\documentclass[a4paper,11pt]{amsart}

\usepackage[english]{babel}
\usepackage{amsmath, enumerate, amsfonts, amssymb, amsthm}

\newtheorem{defin}{\bf Def\mbox{}inition}[subsection]
\newtheorem{theo}[defin]{\bf Theorem}
\newtheorem{prop}[defin]{\bf Proposition}
\newtheorem{lem}[defin]{\bf Lemma}
\newtheorem{cor}[defin]{\bf Corollary}
\newtheorem{clai}[defin]{\bf Claim}
\newtheorem{rem}[defin]{\bf Remark}

\newtheorem{ex}[defin]{\bf Example}
\newtheorem{algo}[defin]{\bf Algorithm}
\newtheorem*{ex*}{\bf Example}
\newtheorem*{clai*}{\bf Claim}
\newtheorem*{rem*}{\bf Remark}
\newtheorem*{nota*}{\bf Notation}
\newtheorem*{prop*}{\bf Proposition}

\newtheorem{theo-intro}{\bf Theorem}
\newtheorem{prop-intro}{\bf Proposition}
\newtheorem{cor-intro}{\bf Corollary}

\newtheorem*{theo:m1*}{\bf Theorem \ref{theo:m1}}

\newcommand{\dps}{\displaystyle}

\newcommand{\C}{\mathbb{C}}
\newcommand{\N}{\mathbb{N}}
\newcommand{\Z}{\mathbb{Z}}

\renewcommand{\k}{\mathbf{k}}
\newcommand{\K}{\mathbf{K}}

\newcommand{\CC}{\mathcal{C}}

\newcommand{\Frac}{\mathrm{Frac}}

\newcommand{\QQ}{\mathcal{Q}}
\newcommand{\modQ}{{\mathrm{mod}\QQ}}
\newcommand{\PP}{\mathcal{P}}
\newcommand{\II}{\mathcal{I}}
\newcommand{\GG}{\mathcal{G}}
\newcommand{\tGG}{\tilde{\mathcal{G}}}
\newcommand{\HHF}{\mathcal{H F}}
\newcommand{\HHP}{\mathcal{H P}}

\newcommand{\spec}{\mathrm{Spec}}

\newcommand{\lm}{\mathrm{lm}} 
\newcommand{\lc}{\mathrm{lc}} 
\newcommand{\lt}{\mathrm{lt}} 
\newcommand{\Supp}{\mathrm{Supp}} 
\newcommand{\Exp}{\mathrm{Exp}} 
\newcommand{\ecart}{\textrm{\'ecart}}

\newcommand{\NF}{\mathrm{NF}} 
\newcommand{\NFMora}{\mathrm{NFMora}}
\newcommand{\PSBmod}{\mathrm{PSBmod}}

\newcommand{\tJ}{\tilde{J}}

\newcommand{\HSf}{\mathrm{HSF}} 
\newcommand{\hHf}{{}^h\mathrm{HF}}   
\newcommand{\aHf}{{}^a\mathrm{HF}} 
\newcommand{\HSp}{\mathrm{HSP}}    
\newcommand{\hHp}{{}^h\mathrm{HP}} 
\newcommand{\aHp}{{}^a\mathrm{HP}}       

\title[Parametric standard basis and degree bound]
{Parametric standard basis,
degree bound and local Hilbert-Samuel function}

\author{Rouchdi BAHLOUL}
\address{Institut Camille Jordan\\
Universit\'e Claude Bernard - Lyon 1\\
43 Boulevard du 11 novembre 1918\\
69622 Villeurbanne\\
France}
\email{bahloul@math.univ-lyon1.fr}

\begin{document}

\begin{abstract}
We propose a general study of standard bases of polynomial
ideals with parameters in the case where the monomial
order is arbitrary. We give an application to the
computation of the stratification by the local Hilbert-Samuel
function. Moreover, we give an explicit upper
bound for the degree of a standard basis for an arbitrary
order and also for the number of the possible affine or
local Hilbert-Samuel functions depending on the
number of variables and the maximal degree of the
given generators.
\end{abstract}


\maketitle


\section*{Introduction and statement of the main results}

In affine algebraic geometry, several (global) objects can be computed
using Gr\"obner bases such as the affine Hilbert polynomial
or free resolutions and parametric Gr\"obner bases may
be seen as a tool for studying these objects under deformations.
In the same way, parametric standard
bases (with respect to local monomial orders) can
be used to study local objects under deformations.

To our knowledge,
most of the existing papers on parametric Gr\"obner or standard
bases concern \emph{global} monomial orders (see
e.g. Lejeune-Jalabert and Philippe \cite{lej-phi},
Gianni \cite{gianni}, Weispfenning \cite{weisp92, weisp03},
Kalkbrenner \cite{kalk}, Montes \cite{montes},
Sato and Suzuki \cite{satosuzuki},
Gonzalez-Vega et al. \cite{gtz}).
In \cite{assi} worked with both global and local orders
(to study flatteners of projections) where the ring
of the coefficients is polynomial.
In \cite{asch2},
Aschenbrenner made a general study of parametric
ideals in power series rings. He also treated the case
where the input generators are polynomials (and
the monomial order is local).
In \cite{jmsj} the author applied parametric standard
bases in rings of differential operators to
study the local Bernstein-Sato polynomial of a deformation
of a hypersurface singularity.

In the present paper, we propose a general study of
parametric standard bases for ideals in some ring
$\CC[x_1, \ldots,x_n]$ where the monomial order on
the $x$-variables is \emph{arbitrary} and the ring
$\CC$ of parameters is also arbitrary.

We shall be concerned both by existencial
and by algorithmic questions. As an application, an algorithm
for computing the stratification by the local Hilbert-Samuel
function is given. Moreover, as an application of a paper by
T. Dub\'e \cite{dube}, we give some bounds for the degree of
standard bases with respect to any monomial order and also for
the number of the possible local or affine Hilbert-Samuel
functions.\\

Before stating the main results, let us introduce some
notations.\\

Throughout the paper,
$\CC$ shall denote an integral domain.
This ring shall be seen as the ring of parameters.
Let $n$ be a positive integer and let $x$ denote the
set $(x_1, \ldots,x_n)$ of indeterminates.\\
Let $\preceq$ be a monomial order on the monomials
$x^\alpha=\prod_i x_i^{\alpha_i}$ ($\alpha \in \N^n$).
We don't suppose $\preceq$ to be a well-ordering (i.e. global).

A \textbf{specialization of} $\CC$ is a ring homomorphism
$\sigma : \CC \to \K$ to some field $\K$. A specialization
$\sigma$ of $\CC$ induces a ring homomorphism $\CC[x] \to \K[x]$
that we shall denote by the same symbol $\sigma$.

The next examples illustrate the situations that we shall
consider in this paper.

\begin{ex}\label{ex:specialization}
\begin{enumerate}
\item
Let $\k \subset \K$ be two fields and $y=(y_1, \ldots, y_m)$ be
a set of indeterminates. For any $y_0\in \K^m$ the map
$(\k[y] \to \K,P \mapsto P(y_0))$ is a specialization
of $\k[y]$ to $\K$. It induces the natural map
$(\k[x,y] \to \K[x], f \mapsto f_{|y=y_0})$.
\item
The previous example is a particular case of the following one.
Given a prime ideal, that is $\PP \in \spec(\K[y])$,
the natural composition map
$\sigma_\PP : \k[y] \to \K[y] \to \K[y]/\PP \subset \Frac(\K[y]/\PP)$
is a specialization. For $y_0=(y_{0,1}, \ldots, y_{0,m}) \in \K^m$,
denote by $m_{y_0}=\sum \K[y](y_j - y_{0,j})$ the maximal ideal
associated with $y_0$. Then $\sigma_{m_{y_0}}$ is identified with
the specialization of (1).
\item
Let $d$ be a positive integer.
Set $q=\binom{n+d}{n}$. This number is the dimension of
the vector space $\oplus_{|\alpha|\le d} \K x^\alpha$ for any field $\K$.
Consider the variables $a=(a_{j,\alpha} | j=1, \ldots, q;
\alpha \in \N^n, |\alpha| \le d)$. Set $N=q^2$. It is
the number of the $a_{j,\alpha}$'s.

For $j=1, \ldots,q$, set
$f_j=\sum_{|\alpha| \le d} a_{j,\alpha} x^\alpha \in \Z[a,x]=\Z[a][x]$.
Let $J=J(n,d)$ denote the ideal of $\Z[a,x]$ generated
by $f_1, \ldots, f_q$.

Put $\CC=\Z[a]$. Let $\K$ be any field. For $a_0 \in \K^N$
we have the natural map $\sigma_{a_0}:\CC \to \K, P(a) \mapsto P(a_0)$.
This specialization induces a map $\sigma_{a_0}:\Z[a,x] \to \K[x]$.

This kind of specialization is interesting because for any
field $\K$ and for any ideal $I$ of $\K[x]$ generated by
polynomials whose degree is at most $d$, there exists
$a_0 \in \K^N$ such that $I=\K[x] \sigma_{a_0}(J)$.
\end{enumerate}
\end{ex}

\subsection{Main results for parametric standard bases}

For a non-zero polynomial $f \in R[x]$ with coefficients in
some ring $R$, $\exp_\preceq(f) \in \N^n$ denotes the leading exponent
of $f$ with respect to $\preceq$, it is defined as the maximum
of the $\alpha$'s such that $x^\alpha$ appears in the development
of $f$.

\begin{theo}\label{theo:main}
Let $J \subset \CC[x]$ and $\QQ \subset \CC$
be finitely generated ideals
such that $J \nsubseteq \CC[x] \cdot \QQ$. 
There exists a finite
set $\GG \subset J$ and finitely many $h_i \in \CC
\smallsetminus \QQ$ such that if we set $h=\prod_i h_i$ then
for any field $\K$ and any specialization $\sigma: \CC \to \K$ such
that $\sigma(\QQ)=\{0\}$ and $\sigma(h)\ne 0$ the following
holds~:
\begin{itemize}
\item
$\sigma(\GG)$ is a $\preceq$-standard basis of $\K[x] \sigma(J)$.
\item
for each $g\in \GG$, $\exp_\preceq(\sigma(g))$ is independent
of $\sigma$.
\end{itemize}
\end{theo}

Notice that the set of the specializations $\sigma$ such that
$\sigma(\QQ)=\{0\}$ and $\sigma(h) \ne 0$ may be empty:
\begin{lem}
Let $\QQ\subset \CC$ be an ideal and $h\in \CC$
then: $(1) \Longrightarrow (2) \iff (3)$, where:
\begin{enumerate}
\item
$h\in \sqrt \QQ$,
\item
For any field $\K$ and any specialization
$\sigma:\CC \to \K$ we have: $\sigma(\QQ)=\{0\} \Rightarrow
\sigma(h)=0$,
\item
$V(\QQ) \subset V(h)$ where $V(\cdot)$ means the affine
scheme defined by (see the notations 0.4).
\end{enumerate}
\end{lem}
\begin{proof}
Assume that $h^i \in \QQ$ for some positive integer $i$.
For a specialization $\sigma : \CC \to \K$, if $\sigma(\QQ)=\{0\}$
then $(\sigma(h))^i=\sigma(h^i)=0$ wich implies $\sigma(h)=0$.
Thus $(1) \Rightarrow (2)$.
Assume $(2)$. Let $\PP \subset \CC$ be prime such that
$\QQ \subset \PP$. Condition (2) applied to $\sigma_\PP$
(as defined in Example~\ref{ex:specialization}(2)) implies
that $\sigma_\PP(h)=0$ which means that $h\in \PP$. Thus we
have $(2) \Rightarrow (3)$.
Conversely assume Condition (3). Let $\sigma$ be a specialization
such that $\sigma(\QQ)=\{0\}$. Then $\ker(\sigma) \in V(\QQ)$.
Therefore $\sigma(h)=0$.
\end{proof}
In 4.3 we shall prove the last implication for $\CC=\k[y]$.

\begin{cor}\label{cor:stratmain}
Assume that $\CC$ is noetherian and
let $J$ be a finitely generated ideal of $\CC[x]$.
There exists a finite set of triples $(\GG_k, \QQ_k, h_k)$
where each $\GG_k \subset J$ is finite,
each $\QQ_k \subset \CC$ is an ideal and each $h_k \in \CC$
and there exists an ideal $\II \subset \CC$ such that
\begin{itemize}
\item
$\spec(\CC)=(\bigcup_k V(\QQ_k) \smallsetminus V(h_k)) \cup V(\II)$,
\item for any specialization $\sigma$ of $\CC$,
if $\sigma(\II)=\{0\}$ then $\sigma(J)=\{0\}$,
\item
for any $k$, for any field $\K$ and any specialization
$\sigma : \CC \to \K$ such that
$\sigma(\QQ_k)=\{0\}$ and $\sigma(h_k)\ne 0$,
\begin{itemize}
\item
$\sigma(J)\ne \{0\}$,
\item
$\sigma(\GG_k)$ is a $\preceq$-standard basis of $\K[x] \sigma(J)$,
\item
for each $g\in \GG_k$, $\exp_\preceq(\sigma(g))$ is independent
of $\sigma$.
\end{itemize}
\end{itemize}
\end{cor}
Here again, $V(\cdot )$ stands for the affine scheme
(see 0.4).

If we form the union of the obtained $\GG_k$ we get a comprehensive
$\preceq$-standard basis $\GG$ (see \cite{weisp92, weisp03, montes}
in the case of a well-ordering $\preceq$):

\begin{cor}\label{cor:compSB}
Let $\CC$ be noetherian and let $J \subset \CC[x]$ be a finitely
generated ideal. There exists a finite set
$\GG \subset J$ such that for any specialization
$\sigma : \CC \to \K$ such that $\sigma(J) \ne \{0\}$,
$\sigma(\GG)$ is a $\preceq$-standard basis of $\K[x] \sigma(J)$.
\end{cor}

\begin{defin}
The ring $\CC$ is called \textbf{detachable} if for any
$h, h_1, \ldots,h_q \in \CC$ there is a finite algorithm for
deciding if $h \in \sum_{j=1}^q \CC \cdot h_j$.
\end{defin}

\begin{prop}\label{prop:detachable}
Suppose that $\CC$ is detachable.
\begin{enumerate}
\item
The set $\GG$ and the elements $h_i$ of Theorem~\ref{theo:main}
can be constructed algorithmically (in a finite number of steps).
\item
Assume that the intersection of two finitely generated ideals
is computable in $\CC$ then the triples $(\GG_k, \QQ_k, h_k)$ and
the ideal $\II$ of Corollary~\ref{cor:stratmain} can be
constructed algorithmically.

Moreover, if for any specialization $\sigma:\CC \to \K$,
$\sigma(J)\ne \{0\}$ then we don't need to assume that the
intersection of ideals in $\CC$ is computable.
\end{enumerate}
\end{prop}

Given a computable field $\k$ and a set of variables
$y=(y_1,\ldots, y_m)$,
then $\k[y]$ and $\Z[y]$ are both detachable.
For $\Z[y]$, see e.g. \cite{ayoub}, \cite{gami} and
\cite{asch-int} (and all the citations in \cite{asch-int}).

\subsection{Constructibility results for Hilbert-Samuel functions}

Let $\k \subset \K$ be two fields where $\k$ is supposed to
be computable.

Let $I$ be an ideal in $\k[x]$.
The affine Hilbert function associated with $I$ is defined as
\[
\N \ni r  \mapsto \aHf_{\k[x]/I}(r)=
\dim_\k(\k[x]_{\le r}/(I \cap \k[x]_{\le r}))
\]
where $\k[x]_{\le r}$ is the vector space
$\oplus_{|\alpha|\le r} \k x^\alpha$.

Given $x_0 \in \K^n$,
let $\K[[x-x_0]]:=\K[[x_1-x_{0,1}, \ldots, x_n-x_{0,n}]]$
denote the ring of formal power series at $x_0$.
The local Hilbert-Samuel function $\HSf_{I,x_0}$ of $I$ at $x_0$
(over $\K$) is defined by:
\[
\N \ni r  \mapsto \HSf_{I, x_0}(r)=
\dim_\K(\K[[x-x_0]]/(\K[[x-x_0]] I + m_{x_0}^{r}))
\]
where $m_{x_0}$ is the maximal ideal of the local ring $\K[[x-x_0]]$.

The notation may seem ambigus if $x_0 \in \k$. In fact for
$x_0 \in \k \subset \K$, the local Hilbert-Samuel function
of $I$ at $x_0$ over $\k$ and the one over $\K$ coincide
(see Lemma~\ref{lem:HS_variable_change}).

There exist numerical polynomials $\aHp_I$ and $\HSp_{I,x_0}$
such that for $r\ge r_0$, $\aHf_I(r)=\aHp_I(r)$ and
$\HSf_{I,x_0}(r)=\HSp_{I,x_0}(r)$ for some $r_0 \in \N$.
These polynomials are called the affine Hilbert polynomial
of $I$ and the local Hilbert-Samuel polynomial of $I$ at $x_0$.

The following is an application of Corollary~\ref{cor:stratmain}.

\begin{cor}\label{cor:stratHS}
There exists an algorithm for computing a finite partition
of $\K^m=\cup W_k$ into constructible sets defined over $\k[x]$
such that for any $W_k$, the map $W_k \ni x_0 \mapsto \HSf_{I,x_0}$ is
constant.
\end{cor}

Let us state another application.
Take the notations of Example~\ref{ex:specialization}(3).

\begin{cor}\label{cor:stratHSgen}
Let $J \subset \Z[x,a]$ be the ideal generated
by the $f_j$'s.
For any field $\K$, there exist a finite partition of $\K^{n+N}$
into constructible subsets $W_k$ with the following properties:
\begin{itemize}
\item
For each stratum $W_k$, and for any $(a_0, x_0)\in W_k$ the local
Hilbert-Samuel function of $J_{|a=a_0} \subset \K[x]$ at $x=x_0$ is constant.
\item
The stratification is defined by ideals in $\Z[a,y]$ that only depend
on the integers $n$ and $d$.
\end{itemize}
\end{cor}

\subsection{Bounds for Standard bases and Hilbert-Samuel functions}

Applying Corollary~\ref{cor:compSB} to
Example~\ref{ex:specialization}(3), one deduces the
existence of a uniform bound $\beta(n,d)$ such that
for any field $\K$ and any ideal in $\K[x]$ generated
by polynomials in $n$ indeterminates of degree at
most $d$, there is a $\preceq$-standard basis whose
elements have degree bounded by $\beta(n,d)$.

In fact, by a direct application of a result by
Dub\'e \cite{dube} (see also the recent generalisation
\cite{al}) we obtain an explicit bound
from which we deduce a bound for the number of the possible
affine or local Hilbert-Samuel functions and polynomials
depending on $n$ and $d$.
This answers some questions by Aschenbrenner in the local
case (see the discussions after Corollary 3.16
and Lemma 3.18 in \cite{asch2}).

Set $D(n,d)=2\bigg( \frac{d^2}{2}+d \bigg)^{2^{n-1}}$.

\begin{prop}\label{prop:boundSB}
Let $d$ and $n$ be positive integers. Let $\preceq$ be any monomial
order on the monomials $x^\alpha=x_1^{\alpha_1}\cdots x_n^{\alpha_n}$.
Given any field $\K$, let $I$ be an ideal of $\K[x_1, \ldots, x_n]$
generated by polynomials of degree at most $d$. Then there exists
a $\preceq$-standard basis of $I$ such that
each element has degree at most $D(n,d)$.
\end{prop}

\begin{prop}\label{prop:boundHS}
Let $d$ and $n$ be positive integers.
There exists a set of functions $\HHF(n,d)$ (from $\N$ to $\N$)
and a set of numerical polynomials $\HHP(n,d)$ that depend
only on $n$ and $d$ such that the following holds.
\begin{itemize}
\item
The cardinality of $\HHP(n,d)$ is
$\dps \binom{nD(n,d)+n}{n}$.
\item
The cardinality of $\HHF(n,d)$ is
\[
\binom{nD(n,d)+n}{n} \cdot
\prod_{k=0}^{nD(n,d)} \left( 1+\binom{k+n-1}{n-1} \right) .
\]
\item
Let $\K$ be a field. Let $I\subset \K[x_1,\ldots,x_n]$ be an ideal
generated by polynomials of degree at most $d$.
\begin{itemize}
\item
$\aHp_{\K[x]/I} \in \HHP(n,d)$ and $\aHf_{\K[x]/I} \in \HHF(n,d)$,
\item
for $x_0 \in \K^n$,
$\HSp_{I,x_0} \in \HHP(n,d)$ and $\HSf_{I,x_0} \in \HHF(n,d)$.
\end{itemize}
\end{itemize}
\end{prop}

\subsection{Main notations}

\begin{itemize}
\item
$\k$: a computable field.
\item
$\K$: an arbitrary field (In many situations we shall
have $\k\subset \K$).
\item
$\langle f_1, \ldots, f_q \rangle$: the ideal generated by
the $f_i$'s.
\item
$x=(x_1,\ldots,x_n)$, $y=(y_1,\ldots,y_m)$: sets of variables.
\item
For $\alpha \in \N^n$ and $\beta\in \N^m$, $x^\alpha:=x_1^{\alpha_1}
\cdots x_n^{\alpha_n}$ and $y^\beta:=y_1^{\beta_1} \cdots y_n^{\beta_n}$.
\item
$\CC$: an integral domain (that may be noetherian or/and detachable).
\item
$\spec(\CC)=\{\PP \subset \CC \ | \ \PP \text{ is a prime ideal}\}$:
the spectrum of $\CC$.
\item
For $S \subset \CC$, $V(S):=\{\PP \in \spec(\CC) \ | \ S\subset\PP\}$:
the affine scheme defined by $S$.
\item
$\sigma$: a specialization to some field $\K$.
\item
For an ideal $J \subset \k[x,y]$ and $y_0 \in \K^m$,
$J_{|y=y_0}:=\K[x] \cdot \{f(x,y_0) \, | \, f(x,y) \in J\}$:
the specialization of $J$ to $y=y_0$.
\item
For an ideal $I \subset \K[x]$ and $x_0 \in \K$,\\
$\aHf_{\K[x]/I}$: the affine Hilbert-function,\\
$\aHf_{\K[x]/I}(r)=\dim_{\K} \K[x]_{\le r}/(I\cap \K[x]_{\le r})$;\\
$\HSf_{I,x_0}$: the local Hilbert-Samuel function at $x=x_0$,\\
$\HSf_{I,x_0}(r)=\dim_{\K} ( \K[[x-x_0]]/
(m_{x_0}^r+\K[[x-x_0]]\cdot I ))$;\\
$\hHf_{\K[x]/I}$: the homogeneous Hilbert function (for $I$ homogeneous),\\
$\hHf_{\K[x]/I}(r)=\dim_{\K} \K[x]_{r}/(I\cap \K[x]_{ r})$
\end{itemize}

\subsection{Structure of the paper}

In section 1, we recall basic facts about standard bases for
polynomial ideals following \cite{singular} and about
Hilbert(-Samuel) functions.
In section 2, we prove the results concerning explicit bounds
(that is Propositions \ref{prop:boundSB} and \ref{prop:boundHS})
since their proof is independent of the rest of the paper.
In section 3, we introduce the notion of pseudo standard basis
modulo some ideal and prove Theorem \ref{theo:main}
and Proposition~\ref{prop:detachable}(1).
In 3.2, we shall propose an alternative method in the
case $\CC=\k[y]$ using usual standard bases.
In section 4, we shall prove Corollary~\ref{cor:stratmain}
and Proposition~\ref{prop:detachable}(2). We shall give
two algorithm that works with a general ring $\CC$ and
one for $\CC=\k[x]$.
In section 5, we shall prove Corollaries \ref{cor:stratHS} and
\ref{cor:stratHSgen}.
We have implemented the algorithm for computing a stratification
with constant local Hilbert-Samuel functions in the computer
algebra system Risa/Asir \cite{asir}.
In section 6, we shall present
some examples computed with our program\footnote{This program
is available on the author's webpage}.

\

\noindent
{\bf Acknowledgements.}
This paper began with some discussions with
Monique Lejeune-Jalabert and Vincent Cossart to whom
I am grateful. I would like to thank Anne Fr\"uhbis-Kr\"uger
for a reading of a the first version of the manuscript
and valuable suggestions.

\section{Recalls on standard bases and Hilbert-Samuel functions}

For 1.1 and 1.2, the reader can refer to chapters 1 and 2
of the book Singular \cite{singular}.

\subsection{Monomial order and associated ring}

As usual, if $\alpha \in \N^n$ then $x^\alpha$ denotes $x_1^{\alpha_1}
\cdots x_n^{\alpha_n}$.
A {\bf monomial order} is a total order $\preceq$ on the
monomials $x^\alpha$ which is compatible with the product,
that is: if $x^\alpha \prec x^{\alpha'}$
then for any $\alpha''$, $x^{\alpha+\alpha''} \prec x^{\alpha'+ \alpha''}$.
An order $\preceq$ is called {\bf global} if $1$ is the
minimal monomial; {\bf local} if $1$ is maximal; mixed otherwise.
In the sequel, we will identify a monomial order with the induced
order on $\N^n$ (which is compatible with the sum).


Let $\k$ be a field.
Let $A$ be a ring with inclusions of rings $\k[x] \subseteq A
\subseteq \k[[x]]$ and let $\preceq$ be a monomial order.
For $f \in A$, write $f=\sum_\alpha c_\alpha x^\alpha$ as a power
series expension. We define the support of $f$ as $\Supp(f)=\{\alpha
\in \N^n | c_\alpha \ne 0\}$. When they make sense, we define the
{\bf leading exponent} of $f$ $\exp(f)=\max_\preceq \Supp(f)$,
the {\bf leading term} $\lt_\preceq(f)=x^{\exp_\preceq(f)}$,
the {\bf leading coefficient} $\lc_\preceq(f)=c_{\exp_\preceq(f)}$ and
the {\bf leading monomial} $\lm_\preceq(f)=\lc_\preceq(f) \lt_\preceq(f)$.
These notions always make sense if $A=\k[x]$. If $A=\k[[x]]$, they
always make sense if $\preceq$ is local.

Now, let us fix a monomial order $\preceq$. Let $R=\k[x]_\preceq$ be
the localization of $\k[x]$ with respect to the multiplicative set
$S_\preceq=\{g \in \k[x] \smallsetminus\{0\} | \exp_\preceq(f)=0\}$.
Notice that $R=\k[x]$ if and only if $\preceq$ is global and
$R=\k[x]_{(0)}$, that is the localization at $0$, if and only if
$\preceq$ is local. In any case we have an inclusion of rings
$\k[x] \subseteq R \subseteq \k[[x]]$. Thus the notations above
apply to the elements of $R$. Notice that if $f\in R$ and
$g \in S_\preceq$ satisfies $gf \in \k[x]$ then
$\exp_\preceq(f)=\exp_\preceq(gf)$.

\subsection{Standard bases in the algebraic situation}

For simplification, we shall forget the subscript $\preceq$.
For the moment $A$ denotes either $\k[x]$ or $R=\k[x]_\preceq$.
Let $J$ be a non zero ideal of $A$. We define the set of leading
exponents $\Exp(J)=\{\exp(f) | f \in J \smallsetminus \{0\}\}$.

\begin{defin}\label{def:Alg_sb}
A finite set $G \subset A$ is called a standard basis of $J$ if
$G\subset J$ and $\Exp(J)=\bigcup_{g \in G} (\exp(g) + \N^n)$.
\end{defin}
By Dickson lemma (see \cite[lemma 1.2.6]{singular}), a standard
basis exists.\\

\noindent
{\bf Remark.}
Assume that $A=\k[x]$.
If $\preceq$ is global we shall use the terminology
{\bf Gr\"obner basis} instead of standard basis.
A Gr\"obner basis generates the ideal but a standard basis
does not in general.

From now on, $A=R=\k[x]_\preceq$.

\begin{defin}[{\cite[Def. 1.6.4]{singular}}]\label{def:Alg_NF}
Let $\mathcal{S}(R)$ denote the set of finite subsets of $R$.
A map $\NF : R \times \mathcal{S}(R) \to R$, $(f, G) \mapsto
\NF(f|G)$ is called a normal form if, for any $f \in R$ and $G
\in \mathcal{S}(R)$, we have
\begin{itemize}
\item[(0)]
$\NF(0 | G)=0$,
\item[(1)]
$\NF(f|G) \ne 0 \Rightarrow \exp(\NF(f|G)) \notin \bigcup_{g\in G}
(\exp(g) +\N^n)$,
\item[(2)]
there exists $u \in R^*=S_\preceq$ and for each $g \in G$, there exists
$a_g \in R$ such that
$r:=u f-\NF(f|G)$ has a standard representation:
$r=\sum_{g \in G} a_g \cdot g$, with
$\exp(r) \succeq \exp(a_g g)$ for all $g$ such that $a_g \ne 0$,
\item[(3)]
if $\{f\} \cup G \subset \k[x]$ then the $a_g$ and $u$ above can
be taken in $\k[x]$. 
\end{itemize}
\end{defin}
\noindent
{\bf Remark.}
This is the definition of a polynomial weak normal form in
the terminology of \cite{singular}.

A normal form always exists: see \cite[1.6, 1.7]{singular} with
NFBuchberger when $\preceq$ is global and NFMora in general.
NFMora is a variant of Mora's division \cite{mora}.

\begin{lem}[{\cite[lemma 1.6.7]{singular}}]
Let $J$ be an ideal of $R$, $G$ be a standard basis of $J$ and
$\NF$ be a normal form.
For any $f \in R$, $f \in J$ if and only if $\NF(f|G)=0$.
\end{lem}
Consequently, $G$ generates $J$ over $R$ (but not over $\k[x]$
in general).

\begin{defin}\label{def:Alg_Spoly}
Let $f,g$ be non zero elements in $R$. Set $\alpha=\exp(f)$,
$\beta=\exp(g)$ and $\gamma=(\gamma_1, \ldots, \gamma_n)$ with
$\gamma_i=\max(\alpha_i, \beta_i)$. We define the $S$-polynomial
(or $S$-function) of $f$ and $g$ as: $S(f,g):=
\lc(g) x^{\gamma-\alpha} f - \lc(f) x^{\gamma-\beta}g$
\end{defin}

\begin{theo}[{\cite[Th. 1.7.3]{singular}}]\label{theo:buch_crit}
Let $J \subset R$ be an ideal and $G$ a finite subset of $J$.
Let $\NF$ be a normal form. The following are equivalent:
\begin{itemize}
\item[(1)]
$G$ is a standard basis of $J$.
\item[(2)]
$\NF(f|G)=0$ for any $f \in J$.
\item[(3)]
Each $f \in J$ has a standard representation with respect to $G$ that is:
there exist some $a_g \in R$ such that
$f=\sum_{g \in G} a_g g$ with $\exp(f) \succeq \exp(a_g g)$ for all
$g$ such that $a_g \ne 0$. 
\item[(4)]
$G$ generates $J$ and for any $g,g' \in G$, $\NF(S(g,g')|G)=0$.
\item[(5)]
$G$ generates $J$ and for any $g_1,g_2 \in G$, there exist some
$a_g \in R$ such that $S(g_1,g_2)=\sum_{g\in G} a_g g$ with
$\exp(S(g_1,g_2)) \succeq \exp(a_g g)$ for all $g$ such that $a_g \ne 0$.
\end{itemize}
\end{theo}
The implications $(4) \Rightarrow (1)$ and $(5) \Rightarrow (1)$ are
usually called {\bf Buchberger's criterion}.
\begin{proof}
The equivalences $(1) \iff \cdots \iff (4)$ are proven in
\cite{singular}. The implication $(3)\wedge(4) \Rightarrow (5)$
is trivial. Let us show that $(5) \Rightarrow (4)$.
Assume by contradiction that for some couple $(g_1,g_2)$,
$\NF((g_1,g_2) | G) \ne 0$. Then by Definition~\ref{def:Alg_NF}(2),
$\exp(\NF(S(g_1,g_2)|G)) \notin \cup_{g\in G} (\exp(g)+\N^n)$.
This contradicts the standard representation in $(5)$.
\end{proof}
The following remark is a direct consequence of Buchberger's criterion.

\begin{rem}\label{rem:ktoK}
Let $\k$ and $\K$ be two fields such that $\k \subseteq \K$.
Let $J$ be an ideal of $\k[x]_\preceq$. If $G$ is a standard basis
of $J$ then it is a standard basis of $\K[x]_\preceq J$.
\end{rem}

\subsection{Hilbert-Samuel function and Standard basis}

Let us start with a remark.
\begin{rem*}
Given an ideal $I$ in $\k[x]$ we defined its local Hilbert-Samuel
function at $0$ as by $\HSf_{I,0}(r)=
\dim_\k(\k[[x]]/(\k[[x]]I + \hat{m}^r))$ where $\hat{m}$ denotes
the maximal ideal of $\k[[x]]$. In the litterature,
one can also find this definition:
$\HSf_{I,0}(r)=\dim_\k(\k[x]_0/(\k[x]_0I + m^r))$ where $\k[x]_0$ is
the localization of $\k[x]$ at $0$ and $m \subset \k[x]_0$
the maximal ideal.
These two definitions coincide.
\end{rem*}
\begin{proof}
Given $r\in \N$, the natural ring homomorphism
$\k[x]_0 \to \k[[x]]/(\k[[x]] I +\hat{m}^r)$ is surjective.
Its kernal is $(\k[[x]](\k[x]_0 I+m^r))\cap \k[x]_0$ and it is
equal to $\k[x]_0 I +m^r$ by faithfull flateness
of $\k[[x]]$ over $\k[x]_0$.
\end{proof}

Given a set $E \subset \N^n$, we define its
{\bf Hilbert-Samuel function} $\HSf_E:\N \to \N$ by
\[
\HSf_E(r)= \mathrm{card}\{\alpha\in \N^n;
\, \alpha\in \N^n \smallsetminus E, \, |\alpha|\le r\}.
\]
A {\bf degree-compatible order} $\preceq$ is a monomial
order such that:
$|\alpha| < |\alpha'| \Rightarrow x^\alpha \prec x^{\alpha'}$
for any $\alpha, \alpha' \in \N^n$. Such an order is global.
A {\bf valuation-compatible order} $\preceq$ is a monomial
order such that:
$|\alpha| > |\alpha'| \Rightarrow x^\alpha \prec x^{\alpha'}$
for any $\alpha, \alpha' \in \N^n$. Such an order is local.

The following is well-known (see e.g.
\cite[Chapt. 9, \S 3, Prop. 4]{clo} and
\cite[Prop. 5.5.7]{singular}).

\begin{lem}\label{lem:rappelsHS}
Let $I$ be an ideal in $\k[x]$.
\begin{enumerate}
\item
If $\preceq$ is a degree-compatible order then
$\aHf_I=\HSf_{\Exp_\preceq(I)}$.
\item
If $\preceq$ is a valuation-compatible order then
$\HSf_{I,0}=\HSf_{\Exp_\preceq(I)}$.
\end{enumerate}
\end{lem}

The next lemma is now trivial (using an affine change
of coordinates, Lemma~\ref{lem:rappelsHS} and
Remark~\ref{rem:ktoK}).

\begin{lem}\label{lem:HS_variable_change}
Let $I$ be an ideal of $\k[x]$ given by generators
$f_1, \ldots, f_q$.
Let $J \subset \k[x,y]$ be the ideal generated
by the $f_i(x+y)$.
Let $\K$ be a field containing $\k$.
Let $x_0 \in \K^n$ and let $\preceq$ be
a valuation-compatible order on the monomials $x^\alpha$.
We have:
\[\HSf_{I,x_0}=\HSf_{\Exp_\preceq(J_{|y=x_0})}.\]
\end{lem}

\section{Bounds for standard bases and Hilbert-Samuel functions}

In this section we shall prove Propositions \ref{prop:boundSB} and
\ref{prop:boundHS}.

\subsection{Bounds for standard bases}

Recall that $\preceq$ is an arbitrary monomial order
on the $x^\alpha$'s.
Let us add a new variable $z$ and consider the following order
$\preceq^z$:

$x^\alpha z^k \prec^z x^{\alpha'} z^{k'} \iff
\begin{cases}
|\alpha|+k < |\alpha'|+k' \text{ or}\\
|\alpha|+k =|\alpha'|+k' \text{ and } x^\alpha \prec x^{\alpha'}
\end{cases}$\\
This order is degree-compatible.

Now we are ready to prove Proposition \ref{prop:boundSB}.
\begin{proof}[Proof of Proposition \ref{prop:boundSB}]
Let $f_1, \ldots, f_q$ be polynomials in $\K[x]$ such that
the degree of $f_j$ is lower than or equal to $d$ for each $j$.
Set $I=\langle f_1, \ldots, f_q\rangle \subset \K[x]$.
Writing $f_j=\sum c_{j,\alpha}x^\alpha$, set
$h(f_j)=\sum c_{j,\alpha}x^\alpha z^{d_j-|\alpha|}$
where $d_j$ denotes the degree of $f_j$.
The following result by Lazard \cite{lazard} is classical
(see, e.g., Exerc. 1.7.6 in \cite{singular}).
\begin{lem}\label{lem:classical}
Let $G$ be a homogeneous $\preceq^z$-standard basis of the
homogeneous ideal $\K[x,z]\{h(f_1), \ldots, h(f_q)\}$.
Then $G_{|z=1}$ is a $\preceq$-standard basis of $I$.
\end{lem}
The ideal $I'=\langle h(f_1), \ldots, h(f_q)\rangle$ is
a homogeneous ideal of $\K[x,z]$ generated by homogeneous
polynomials of degree bounded by $d$.
Applying \cite[Theorem 8.2]{dube} by T. Dub\'e, one may choose
$G$ such that the degree of its elements is bounded
by $D(n,d):=2((d^2/2)+d)^{2^{n-1}}$. Therefore the elements
of $G_{|z=1}$ have their degree bounded by $D(n,d)$.
\end{proof}

\subsection{Bounds for affine and local Hilbert-Samuel functions}

Let us begin with a basic combinatorial result.
\begin{lem}\label{lem:combin}
Let $\delta$ be in $\N$. The cardinality of the
following set is $\binom{n+\delta}{n}$:
\[
\{(b_1, \ldots, b_n)\in\N^n\} | \
b_n\le b_{n-1}\le \cdots \le b_1\le \delta\}.
\]
\end{lem}
\begin{proof}
For $(b_1, \ldots, b_n)$ in this set one can associate the
following monomial: $x_1^{b_1-b_2} \cdots x_{n-1}^{b_{n-1}-b_n} x_n^{b_n}$.
This induces a bijective map from our set to the set of
the monomials $m$ of degree $\deg(m)\le \delta$.
It is well-know that the cardinality of the latter
is $\binom{n+\delta}{n}$ (see e.g. Lemma 4 page 438
in \cite{clo}).
\end{proof}

Here is another technical lemma.
\begin{lem}\label{lem:borne}
Let $J$ be a monomial ideal in $\K[x_1, \ldots,x_n]$.
Let $G$ be a finite
set of monomials generating $J$ and let
$\delta=\max\{\deg(m) | m\in G\}$. Then for any
$r \ge n\delta$, $\aHf_{\K[x]/J}(r)=\aHp_{\K[x]/J}(r)$
\end{lem}
\begin{proof}
Set $G=\{m_1, \ldots, m_q\}$.
For $t=(t_1, \ldots t_k)$ such
that $1\le k \le q$ and $1 \le t_1 < \cdots < t_k \le q$.
let $M_t$ be the ideal
generated by $m_{t_1}, \ldots, m_{t_k}$. Finally, for $r\in \N$,
set $M_{t,r}=\{x^\alpha ; |\alpha| \le r, x^\alpha \in M_t\}$.
Applying the inclusion-exclusion principle, one obtains
\begin{eqnarray*}
\aHf_{\K[x]/J}(r)
&=& \binom{r+n}{n}-\mathrm{card} (M_{1,r}\cup \cdots \cup M_{q,r})\\
&=& \binom{r+n}{n}-\sum_{k=1}^q (-1)^{k-1}
    \sum_{1\le t_1 < \cdots < t_k \le q} \mathrm{card} (M_{(t_1,\ldots,t_k),r}).
\end{eqnarray*}
Since $M_{(t_1,\ldots,t_k)}$ is generated by
$\mathrm{lcm}\{m_{t_1},\ldots,m_{t_k}\}$, we have
$\mathrm{card}(M_{(t_1,\ldots,t_k),r})=\binom{r+n-e}{n}$ for
every $r\ge e$ where $e$ is the degree of this common multiple.
Since $n\delta$ is a bound for the degree of all the common multiples,
we conclude that $\aHf_{\K[x]/J}(r)$ is polynomial for
$r\ge n\delta$.
\end{proof}

Let us recall some facts from Dub\'e's paper \cite{dube}.
For this, we recall that given a homogeneous ideal $J$
in $\K[x_1,\ldots, x_n]$ one can define the (homogeneous)
Hilbert function
\[
\hHf_{\K[x]/J}(r)=\dim_{\K}(\K[x]_r/\K[x]_r \cap J)
\]
where $\K[x]_r=\oplus_{|\alpha|=r} \K x^\alpha$.
Knowing the affine Hilbert function or the homogeneous one
is equivalent since we have~:
$\hHf_{\K[x]/I}(r)=\aHf_{\K[x]/I}(r)-\aHf_{\K[x]/I}(r-1)$
and $\aHf_{\K[x]/I}(r)=\sum_{k=0}^r \hHf_{\K[x]/I}(k)$.
Dub\'e defines the Macaulay constants $(b_0, \ldots, b_n) \in \N^n$
for any homogeneous ideal $J$. These numbers are uniquely determined
and they have some properties~:
\begin{itemize}
\item
$b_0 \ge b_1 \ge \cdots \ge b_n$.
\item
$b_0$ is equal to $\min \{b\in \N | \forall r\ge b,
\hHf_{\K[x]/J}(r)=\hHp_{\K[x]/J}(r)\}$.
\item
$\dps \hHp_{\K[x]/I}(r)=\binom{r+n}{n}-1-
\sum_{k=1}^n\binom{r-b_k+k-1}{k}$.
\end{itemize}
This shows that the constants $b_1, \ldots,b_n$ uniquely
determines and are uniquely determined by the (homogeneous)
Hilbert polynomial.

Now let us prove Proposition~\ref{prop:boundHS}.
\begin{proof}[Proof of Prop. \ref{prop:boundHS}]
We shall begin by the local case.
Recall that we start with an ideal $I \subset \K[x]$
that admits a finite set of generators whose degree
is bounded by $d$. 
By an affine change of coordinates, we are
reduced to the case where $x_0=0$.
Let us consider a valuation-compatible order.

By Proposition~\ref{prop:boundSB}, $I$ admits a standard basis
$G$ such that the degree of each element is bounded by $D:=D(n,d)$.
Let $M$ be the monomial ideal generated by the leading monomials
of the $g$'s in $G$. By Lemma~\ref{lem:rappelsHS},
$\HSf_{\K[x]/I}=\aHf_{\K[x]/M}$.
By Lemma~\ref{lem:borne}, for $r\ge nD$,
$\aHf_{\K[x]/M}(r)=\aHp_{\K[x]/M}(r)$.

Using the recalls of Dub\'e's results,
we have that $\hHf_{\K[x]/M}$ and then $\aHf_{\K[x]/M}$
is uniquely determined by some tuple $(b_1, \ldots, b_n)$ such that
$nD \ge b_1 \ge \cdots \ge b_n$. By Lemma~\ref{lem:combin},
the number of these tuples is $\binom{n+nD}{n}$. This proves
the part concerning the local Hilbert-Samuel polynomial.
Now again by Dub\'e's results, $b_0=\min \{b\in \N | \forall r\ge b,
\aHf_{\K[x]/M}(r)=\aHp_{\K[x]/M}(r)\}$.

Since $b_0 \le nD$, $\aHf_{\K[x]/M}(r)$ is determined
for all $r\ge nD$.
It remains to count the number of possible values that may
be taken by $\aHf_{\K[x]/M}(r)$ for $0 \le r<nD$.
For a given $0 \le r<nD$, $\aHf_{\K[x]/M}(r)$ may
be in $\{0, \ldots, \binom{r+n-1}{n-1}\}$. Therefore
it may take
$\binom{r+n-1}{n-1}+1$ possible values.
Taking the product for all $0\le r < nD$ we obtain the
bound for the number of the possible
local Hilbert-Samuel functions.

Now, the proof concerning the affine Hilbert function and polynomial
is the same providing the use of a degree-compatible order
instead of a valuation-compatible one.
\end{proof}

\section{Parametric standard bases}

In this section we shall be concerned by the proof of
Theorem~\ref{theo:main} and Proposition~\ref{prop:detachable}(1).
In 3.1 we shall treat the case of a general $\CC$ and
in 3.2 we shall propose another method when $\CC=\k[y]$.

\subsection{General case: an analogue of pseudo standard bases}

Recall that $\CC$ is an integral domain
and $\preceq$ is a monomial order on the $x^\alpha$'s.
For $f\in \CC[x]\smallsetminus\{0\}$,
we can define its leading exponent $\exp_\preceq(f) \in \N^n$,
its leading term $\lt_\preceq(f)=x^{\exp_\preceq(f)}$,
its leading coefficient $\lc_\preceq(f) \in \CC$ and
leading monomial $\lm_\preceq(f)=\lc_\preceq(f) \cdot \lt_\preceq(f)$.
In the sequel we shall forget the subscript $\preceq$ and write
$\exp(f)$ for $\exp_\preceq(f)$.

Set $S_\preceq=\{f \in \CC[x] | \exp(f)=0 \textrm{ and }
\lc(f)=1\}$ then define $R=S_\preceq^{-1} \CC[x]$ as the localization
w.r.t. $S_\preceq$.

\begin{defin}[See \cite{gp} or {\cite[pages 124-125]{singular}}]
\
\begin{itemize}
\item
As in Def. \ref{def:Alg_NF}, $\mathrm{S}(R)$ denotes the set of
finite subsets of $R$.
A map $\NF : R \times \mathcal{S}(R) \to R$, $(f, G) \mapsto
\NF(f|G)$ is called a pseudo normal form if, for any $f \in R$
and $G \in \mathcal{S}(R)$, we have
\begin{itemize}
\item[(0)]
$\NF(0 | G)=0$
\item[(1)]
$\NF(f|G) \ne 0 \Rightarrow \exp(\NF(f|G)) \notin
\bigcup_{g\in G} (\exp(g) +\N^n)$
\item[(2)]
There exists $u \in R$ such that $\lm(u)$
is of the form $\lm(u)=\prod_{g\in G} \lc(g)^{d_g} \cdot x^0$
with $d_g \in \N$,
and for each $g \in G$, there exists
$a_g \in R$ such that
$r:=u f-\NF(f|G)$ has a standard representation:
$r=\sum_{g \in G} a_g \cdot g$, with
$\exp(r) \succeq \exp(a_g g)$ for all $g$ such that $a_g \ne 0$.
\item[(3)]
If $\{f\} \cup G \subset \k[x]$ then the $a_g$ and $u$ above can
be taken in $\CC[x]$. 
\end{itemize}
\item
Given a non-zero ideal $J \subset R$, a pseudo standard
basis is a finite set $G \subset J$ satisfying
$\Exp(J)=\bigcup_{g\in \GG} (\exp(g) + \N^n)$.
\end{itemize}
\end{defin}
Notice that our definition of a pseudo standard basis is
slightly different to the one given in \cite{singular, gp}.

Pseudo normal forms exist (NFMora in \cite{singular} is one)
and pseudo standard bases also (by Dickson lemma).



Now let us generalize these constructions.
In the sequel $\QQ \subset \CC$ is a given ideal, not
necessarily prime.
Given $f \in R=\CC[x]_\preceq$, we define $\exp^\modQ(f):=
\exp(f \modQ)$, where $f \modQ$ means the class of $f$ in
$\CC/\QQ[x]_\preceq$ viewed in $(\CC/\QQ)[[x]]$.
We define $\lt^\modQ(f):=x^{\exp^\modQ(f)}$.
Then $\lc^\modQ(f)$ denotes the coefficient (in $\CC$) of
$\lt^\modQ(f)$ in the expension of $f$, finally $\lm^\modQ(f):=
\lc^\modQ(f) \lt^\modQ(f)$.

Now for an ideal $J \subset R$ such that $J \nsubseteq R\QQ$,
we define
$\Exp^\modQ(J)=\{\exp^\modQ(f) | f \in J \smallsetminus R \QQ\}$.

\begin{defin}\label{def:pseudoSBmodQ}
A pseudo standard basis of $J$ modulo $\QQ$ is a finite set
$\GG \subset J$ such that $\Exp^\modQ(J)=\bigcup_{g \in \GG}
(\exp^\modQ(g) + \N^n)$.
\end{defin}
\begin{rem*}
Such a set exists by Dickson lemma again.
Notice that if $\QQ=(0)$, we recover the notion of a
pseudo standard basis.
\end{rem*}

\begin{defin}\label{def:pnfmodQ}
A pseudo normal form $\NF(\cdot |_\QQ \cdot )$
modulo $\QQ$ is a map
$\NF(\cdot |_\QQ \cdot ) : R \times \mathcal{S}(R) \to R$,
$(f, G) \mapsto \NF(f|_\QQ G)$
such that for any $f \in R$ and $G \in \mathcal{S}(R)$,
we have
\begin{itemize}
\item[(0)]
$\NF(q |_\QQ G) \in R \QQ$ for all $q\in R\QQ$
\item[(1)]
$\NF(f|_\QQ G) \notin R\QQ \Rightarrow
\exp^\modQ(\NF(f|_\QQ G)) \notin
\bigcup_{g\in G} (\exp^\modQ(g) +\N^n)$
\item[(2)]
There exist some $a_g \in R$, $q \in R \QQ$, and $u \in R$ such that
$\lm^\modQ(u)=\prod_{g\in G} (\lc^\modQ(g))^{d_g} \cdot x^0$
with $d_g \in \N$ and

$r:=u f- \NF(f |_\QQ G) = \sum_{g\in G} a_g g + q$\\
with
$\exp^\modQ(r) \succeq \exp^\modQ(a_g g)$ for
all $g$ such that $a_g \ne 0$.
\item[(3)]
If $\{f\} \cup G \subset \CC[x]$ then the $a_g$ and $u$ and $q$
above can be taken in $\CC[x]$. 
\end{itemize}
\end{defin}

We define the $S$-function modulo $\QQ$:
\begin{defin}\label{def:Alg_SpolymodQ}
Let $f,g \in R \smallsetminus R\QQ$. Set $\alpha=\exp^\modQ(f)$,
$\beta=\exp^\modQ(g)$ and $\gamma=(\gamma_1, \ldots, \gamma_n)$ with
$\gamma_i=\max(\alpha_i, \beta_i)$. We define the $S$-polynomial
(or $S$-function) of $f$ and $g$ modulo $\QQ$ as:
$S^\modQ(f,g):= \lc^\modQ(g) x^{\gamma-\alpha} f -
\lc^\modQ(f) x^{\gamma-\beta}g$
\end{defin}

As for standard bases, we have a characterization of pseudo
standard bases in terms of pseudo normal forms and $S$-polynomials.

\begin{prop} \label{prop:pseudo_buch_crit}
Let $J \subset R$ be an ideal and $G$ a finite subset of $J$.
Let $\NF(\cdot |_\QQ \cdot )$ be a pseudo normal
form modulo $\QQ$.
The following are equivalent:
\begin{itemize}
\item[(1)]
$G$ is a pseudo standard basis of $J$ modulo $\QQ$.
\item[(2)]
$\NF(f|_\QQ G) \in R \QQ$ for any $f \in J$.
\item[(3)]
For any $f \in J$, there exists $a_g\in R$ for all $g \in G$,
$q \in R\QQ$, and $u \in R$ with $\lm^\modQ(u)$ being a product
of $\lc^\modQ(g)$ ($g \in G$) such that:
$u f=\sum_{g \in G} a_g g +q$ with $\exp^\modQ(f) \succeq
\exp^\modQ(a_g g)$ for all $g$ such that $a_g \ne 0$. 
\item[(4)]
For any $f\in J$, there exists $u$ as above such that
$uf \in RG+R\QQ$ and
for any $g,g' \in G$, $\NF(S^\modQ(g,g')|_\QQ G) \in R \QQ$.
\end{itemize}
\end{prop}
\begin{proof}
Let us prove $(1) \Rightarrow (2)$.
Assume $(1)$ and
by contradiction let $f \in J$ be such that
$\NF(f|_\QQ G) \notin R\QQ$.
Then by Definition~\ref{def:pnfmodQ}(1),
$\exp^\modQ(\NF(f|_\QQ G)) \notin
\bigcup_{g\in G} (\exp^\modQ(g) +\N^n)$.
By Definition~\ref{def:pnfmodQ}(2), $\NF(f|_\QQ G) \in J + R\QQ$
therefore $\exp^\modQ(\NF(f|_\QQ G)) \in \Exp^\modQ(J)$.
But this contradicts $(1)$.
The proof of $(2) \Rightarrow (3) \Rightarrow (1)$
is a direct application of the definitions.
Moreover
$(3)$ implies the first part of $(4)$ and $(2)$ implies the second
one. It remains to prove (for example) $(4) \Rightarrow (1)$.
For this, let us introduce some extra notations.
For $f \in R$ let us denote by $(f)_\QQ$ its image
by the natural map $\CC[x]_\preceq \to \CC/\QQ[x]_\preceq \to
\Frac(\CC/\QQ)[x]_\preceq$. Let us denote by $(J)_\QQ$ the
ideal generated by $\{(f)_\QQ | \ f\in J\}$.

Now, let us assume $(4)$ and by contradiction suppose
that $(1)$ is not true. There exists $f\in J \smallsetminus R\QQ$
such that $\exp^\modQ(f) \notin \cup_{g\in G} (\exp^\modQ(g)+\N^n)$.
This implies that $\exp((f)_\QQ) \notin \cup_{g\in G}
(\exp(g)_\QQ +\N^n )$. Hence $(G)_\QQ$ is not a standard basis
of $(J)_\QQ$.

The first part of $(4)$ implies that $(G)_\QQ$ generates $(J)_\QQ$.
The second part of $(4)$ combined with
Definition~\ref{def:pnfmodQ}(2) implies that for all $g,g' \in G$,
\[S((g)_\QQ, (g')_\QQ)=(S^\modQ(g,g'))_\QQ = \sum_g (a_g)_\QQ \cdot (g)_\QQ\]
with $\exp(S((g)_\QQ, (g')_\QQ)) \succeq \exp((a_g)_\QQ \cdot (g)_\QQ)$,
hence by Theorem~\ref{theo:buch_crit}(5),
$(G)_\QQ$ is a standard basis of $(J)_\QQ$.
Contradiction.
\end{proof}

Given $f \in \CC[x]$, we define the \'ecart modulo $\QQ$:
$\ecart^\modQ(f):=\ecart(f \modQ)=\deg(f \modQ)- \deg \lt(f \modQ)$.

\begin{algo}[$\NFMora^{\mod \bullet}(\bullet | \bullet)$]
\label{algo:NFMora}
\rm
\

\noindent
{\rm \bf Input:} $f \in \CC[x]$, $G \subset \CC[x]$: a finite set,
$\QQ \subset \CC$: an ideal.\\
{\rm \bf Output:} $h=\NFMora^\modQ(f | G) \in \CC[x]$ a pseudo normal form
of $f$ w.r.t. $G$ modulo $\QQ$.\\

\noindent
$\bullet$ $h:=f$;\\
$\bullet$ $T:=G$;\\
$\bullet$ While $(h \notin \CC[x] \QQ \textrm{ and }
T_h:=\{g \in T \text{ such that } \lt^\modQ(g)|\lt^\modQ(h)\} \ne
\emptyset)$\\
\hspace*{.7cm} $\bullet$ Choose $g \in T_h$ with $\ecart^\modQ(g)$ minimal;\\
\hspace*{.7cm} $\bullet$ If $(\ecart^\modQ(g) > \ecart^\modQ(h))$
then $T:=T \cup\{h\}$;\\
\hspace*{.7cm} $\bullet$ $h:=S^\modQ(h,g)$;\\
$\bullet$ Return $h$.
\end{algo}

All the definitions were made in order to have the following
equality:
$\NFMora^\modQ(f| G) \modQ=\NFMora(f \modQ | G \modQ)$.
This proves both termination and correctness of this algorithm.
Moreover this equality proves that
$\NFMora^{\mod \QQ}(\bullet | \bullet)$ is pseudo normal
form modulo $\QQ$.

To be complete, let us give a generalisation of the
algorithm ``Standard''
(see Algorithm 1.7.1 in \cite{singular}).
\begin{algo}[$\mathrm{Standard}^{\mod \bullet}(\bullet, \bullet)$]
\label{algo:standard}
\rm
\

\noindent
{\rm \bf Input:} $G \subset R$: a finite set,
$\QQ \subset \CC$: an ideal,
$\NF$: a pseudo normal form modulo $\QQ$.\\
{\rm \bf Output: } $S:=\mathrm{Standard}^\modQ(G,\NF)$
a pseudo standard basis
of the ideal $R G$ modulo $\QQ$.

\noindent
$\bullet$ $S:=G$;\\
$\bullet$ $P:=\{(f,g)| f,g \in S, f\ne g\}$;\\
$\bullet$ While $P\ne \emptyset$\\
\hspace*{.7cm} $\bullet$ Choose $(f,g)\in P$;\\
\hspace*{.7cm} $\bullet$ $P:=P\smallsetminus \{(f,g)\}$;\\
\hspace*{.7cm} $\bullet$ $h:=\NF(S^\modQ(f,g) |_\QQ S)$;\\
\hspace*{.7cm} $\bullet$ If $(h\notin R\QQ)$ then
($P:=P\cup \{(h,f) | f\in S\}$; $S:=S\cup \{h\}$);\\
$\bullet$ Return $S$.
\end{algo}

\begin{clai}\label{claim:final}
Let $G$ be a finite system of generators of $J$.
The set $\GG=\mathrm{Standard}^\modQ(G,\NFMora^\modQ)$
is a pseudo standard basis of $J$ modulo $\QQ$. 
\end{clai}
\begin{proof}
First, notice that $\mathrm{Standard}^\modQ(G,\NFMora^\modQ)$
terminates because,
$\mathrm{Standard}^\modQ(G,\NFMora^\modQ) \modQ=
\mathrm{Standard}(G\modQ, \NFMora)$.
Now let us prove the proposition. The algorithm
$\mathrm{Standard}^\modQ$ terminates when the set $P$ of 
pairs is empty. This set becomes empty when
$\NF(S^\modQ(f,g) |_\QQ S)$ is in $R\QQ$ for all $(f,g) \in P$.
Thus the output $\GG$ satisfies
Condition (4) of Proposition~\ref{prop:pseudo_buch_crit}.
Thus, $\GG$ is a pseudo standard basis of $J$
modulo $\QQ$.
\end{proof}

\begin{prop}\label{prop:final}\
\begin{enumerate}
\item
Let $\GG$ be a pseudo standard basis of $J$ modulo $\QQ$
and let $h=\prod_{g \in \GG} \lc^\modQ(g)$.
For any field $\K$ and any specialization $\sigma :\CC \to \K$
such that $\sigma(\QQ)=\{0\}$ and $\sigma(h) \ne 0$:
\begin{itemize}
\item
$\sigma(\GG)$ is a $\preceq$-standard basis of $\K[x]\sigma(J)$.
\item
for each $g\in \GG$, $\exp(\sigma(g))=\exp^\modQ(g)$.
\end{itemize}
\item
Moreover if $J$ is generated by a set $G \subset \CC[x]$
then it is possible to construct $\GG$ inside $\CC[x] G$.
\end{enumerate}
\end{prop}

Proving this proposition proves Theorem~\ref{theo:main}
\begin{proof}
\begin{enumerate}
\item
Recall that $J$ is an ideal of $\CC[x]$ generated by a given
finite set $G$, and $\QQ \subset \CC$ is
an ideal such that $J \nsubseteq \CC[x]\QQ$.
Fix a field $\K$.

Let $\GG$ be a pseudo standard basis of $J$ modulo $\QQ$.
Let $h \in \CC$ be the product of the $\lc^\modQ(g)$ with $g\in \GG$.
Let $\Sigma$ be the set of the specializations $\sigma:\CC \to \K$
such that $\sigma(\QQ)=\{0\}$ and $\sigma(h)\ne 0$.
For $\sigma \in \Sigma$ and for any $g \in \GG$,
$\sigma(\lc^\modQ(g))\ne 0$ and $\exp^\modQ(g)=\exp(\sigma(g))$
which proves the constancy of $\exp(\sigma(g))$ over
$\sigma \in \Sigma$.

Take $\sigma \in \Sigma$.
Following
$\NFMora^\modQ$ and $\NFMora$ step by step,
we obtain that $\sigma(\NFMora^\modQ(S^\modQ(g,g') | \GG))$
is equal to $\NFMora(S(\sigma(g), \sigma(g')) | \sigma(\GG))$
and it is $0$ for all $g,g' \in \GG$ 
by Prop.~\ref{prop:pseudo_buch_crit}(4).

Proposition~\ref{prop:pseudo_buch_crit}(4)
implies that $\sigma(\GG)$ generates $\K[x]\sigma(J)$.
Thus
Buchberger's criterion (Theorem~\ref{theo:buch_crit})
implies
that $\sigma(\GG)$ is a standard basis of $\K[x]\sigma(J)$.
\item
By definition $S^\modQ(f,g) \in \CC[x]f + \CC[x]g$. Moreover,
a pseudo normal form modulo $\QQ$ $\NF(\cdot |_\QQ \cdot)$
(see Condition (4) in
Definition~\ref{def:pnfmodQ}) outputs an element that
is a combination over $\CC[x]$ of the inputs (it is
obviously true for $\NFMora^\modQ$). Finally, in the
algorithm $\mathrm{Standard}^\modQ$, if the inputs are
in $\CC[x]$ then so are the outputs.
\end{enumerate}
\end{proof}

In order to conclude this part, it remains to prove
Proposition~\ref{prop:detachable}(1).
\begin{proof}[Proof of Prop.~\ref{prop:detachable}(1)]
Suppose that $\CC$ is detachable.
Then it is clear that all the ``objects'' modulo $\QQ$
can be computed (such as $\exp^\modQ(f)$, $S^\modQ(f,g)$).
Thus, given $G \subset \CC[x]$, the set
$\GG=\mathrm{Standard}^\modQ(G,\NFMora^\modQ)$ can
be computed in a finite number of steps.
\end{proof}

\begin{rem}\label{rem:modQ}
\rm
Suppose that $\CC=\k[y]$ with a computable field $\k$.

From an algorithmic point of view
all the ``objects'' $\modQ$ (such as $\exp^\modQ$)
can be computed in the following way.
We consider a monomial order on the
$y^\beta$, say $\le_0$ and compute a standard or Gr\"obner basis
of $\QQ$, say $G_0$. Then we consider a monomial order,
say $\le$, on the monomials $x^\alpha y^\beta$ whose restriction
to $y^\beta$ is $\le_0$ (for example the block order
$(\preceq, \le_0)$). Then e.g. by Buchberger's criterion,
$G_0$ is a standard basis
of $\k[x,y] \QQ$ w.r.t. $\le$. Given $f \in \k[x,y]$, we
compute a normal form $r=\NF_\le(f|G_0)$ and we get
$\exp_\preceq^\modQ(f)=\exp_\preceq(r)$.
\end{rem}

Summing up the results above we get the following
algorithm (when $\CC$ is detachable).

\begin{algo}[$\mathrm{PSBmod}$]
\label{algo:psbmod}
\rm
\

\noindent
{\rm \bf Input:} $G \subset \CC[x]$: a finite set,
$\QQ \subset \CC$: an ideal.\\
{\rm \bf Output:} $\PSBmod(G,\QQ)=(\GG, H)$ where
$\GG$ is a pseudo standard basis of $\langle G \rangle$
modulo $\QQ$,
$H \subset \CC\smallsetminus \QQ$ is a finite set.

\noindent
$\bullet$ $H:=\emptyset$; $\GG:=\emptyset$;\\
$\bullet$ if $G \subset \CC[x]\QQ$ then Return $(\GG ,H)$;\\
$\bullet$ $\GG:=\mathrm{Standard}^\modQ(G,\NFMora^\modQ)$;\\
$\bullet$ for $g\in \GG$ do
 $(H:=H \cup \{\lc^\modQ(g)\})$;\\
$\bullet$ Return $(\GG, H)$.
\end{algo}

If $G \subset \CC[x]\QQ$ then the output is $(\emptyset, \emptyset)$,
otherwise we get $(\GG, H)$ and setting $h$ as
the product of the elements of $H$, we have that $\GG$
specializes to a standard basis for all $\sigma :\CC \to \K$
such that $\sigma(\QQ)=\{0\}$ and $\sigma(h) \ne 0$.

\subsection{The case $\CC=\k[y]$: with standard bases}

In this paragraph, we give an alternative method for
computing pseudo standard bases modulo some $\QQ$ in
the particular case where $\CC=\k[y]$.

Denote by $\tJ=J+\k[x,y] \QQ$. Let $\le_0$ be a monomial order on
the $y^\beta$. For simplicity, we assume that $\le_0$ is \emph{global}.
(In fact, things work even if $\le_0$ is not global: the proof
of Proposition~\ref{prop:main2} would need a slight modification.)

Define a block order on $x^\alpha y^\beta$ as
$\le=(\preceq, \le_0)$, here $\preceq$ is the monomial order
on $x^\alpha$ used from the beginning.

\noindent
{\bf Note.} For an element $f \in \k[x,y]$, we will work with two
types of leading exponents (and of leading terms, coefficients,
etc): $\exp_{\preceq}(f) \in \N^{n}$ and $\exp_{\le}(f)\in \N^{n+m}$.

\begin{rem}\label{rem:exp_le}
For any $f\in \k[x,y]$, $\exp_\le(f)=(\exp_\preceq(f),
\exp_{\le_0}(\lc_\preceq(f)))$.
\end{rem}

\

Let $G$ be a standard basis of $\tilde{J}=J+ \k[x,y] \cdot \QQ$
w.r.t. $\le$.

\begin{prop} \label{prop:main2}
The set $\tGG=G \smallsetminus \k[x,y]\cdot \QQ$ is a
$\preceq$-pseudo standard basis of $\tJ$ modulo $\QQ$.
\end{prop}

\begin{proof}

Take $f \in \tJ$ such that $f \notin \k[x,y] \QQ$.
We are going to prove that $\exp_\preceq^\modQ(f)\in
\exp_\preceq^\modQ(g)+\N^n$ for some
$g\in G \smallsetminus \k[x,y]\QQ$.
Since $\QQ \subset \tJ$, we may
assume $\lc^\modQ_\preceq(f)=\lc_\preceq(f)$.

Let $c \in \k[y]$ be a normal form of $\lc_\preceq(f)$ with
respect to a $\le_0$-Gr\"obner basis of $\QQ$.
Since $\lc_\preceq(f) - c$ is in
$\QQ \subset \tJ$, we may assume that
\[(\star) \qquad
\exp_{\le_0}(\lc_\preceq(f)) \notin \Exp_{\le_0}(\QQ).\]
By definition of $G$, there exists $g \in G$ such that
$\exp_\le(f) \in \exp_\le(g) +\N^{n+m}$.
By Remark~\ref{rem:exp_le}, this implies
$\exp_{\le_0}(\lc_\preceq(f)) \in \exp_{\le_0}(\lc_\preceq(g)) + \N^m$.
Relation $(\star)$ implies $\lc_\preceq(g) \notin \QQ$,
i.e. $\lc_\preceq(g)=\lc_\preceq^\modQ(g)$. Therefore
$g\notin \k[x,y]\QQ$.
By Remark~\ref{rem:exp_le} again, we have
$\exp_\preceq(f) \in \exp_\preceq^\modQ(g) +\N^n$.
\end{proof}

Now let us define $\GG \subset J$ as follows.
For each element $\tilde{g} \in \tGG$ let $g \in J$
be such that $g-\tilde{g} \in \k[x,y]\cdot \QQ$.
We define $\GG$ as the set of these $g$ for $\tilde{g} \in \tGG$.
The set $\GG$ is not uniquely determined of course.

As a trivial consequence of the definition of $\tJ$ we obtain:

\begin{cor}\label{cor:main2}
$\GG$ is a $\preceq$-pseudo standard basis of $J$
modulo $\QQ$.
\end{cor}

Hence, this ends the second proof of Theorem~\ref{theo:main}.

Now in order to end this part, we have to propose
an algorithmic construction for such a $\GG$.
We think that the smplest way is to construct in parallel
the sets $\tGG$ and $\GG$. For this, we propose
a modification of the algorithm Standard.

\begin{algo}[$\mathrm{ModifiedStandard}(\bullet, \bullet, \bullet)$]
\label{algo:standardbis}
\rm
\

\noindent
{\rm \bf Input:} $G_1, G_2 \subset \k[x,y]$: finite sets,
$\NF$: a normal form.\\
{\rm \bf Output: } $S:=\mathrm{ModifiedStandard}(G_1,G_2, \NF)$
where $S=\{(\tilde{g}_1, g_1), \ldots, (\tilde{g}_s, g_s)\}$
is a finite set such that $\{\tilde{g_1}, \ldots, \tilde{g}_s\}$
is a standard basis of $\langle G_1 \cup G_2 \rangle$ and for all $i$,
$g_i \in \langle G_1 \rangle$ and $\tilde{g}_i -g_i \in
\langle G_2 \rangle$.

\noindent
$\bullet$ $S:=\bigcup_{g\in G_1} \{(g, g)\} \cup
\bigcup_{g\in G_2}\{(g, 0)\}$;\\
$\bullet$ $P:=\{((\tilde{f}, f),(\tilde{g},g))| \
(\tilde{f}, f),(\tilde{g},g) \in S, \tilde{f}\ne \tilde{g}\}$;\\
$\bullet$ While $(P\ne \emptyset)$\\
\hspace*{.5cm} $\bullet$ choose $((\tilde{f}, f) ,(\tilde{g}, g))\in P$;\\
\hspace*{.5cm} $\bullet$ $P:=P\smallsetminus \{ ((\tilde{f}, f) ,(\tilde{g}, g)) \}$;\\
\hspace*{.5cm} $\bullet$ $\tilde{h}:=\NF(S(f,g) | S)$;\\
\hspace*{.5cm} $\bullet$ If $(\tilde{h} \ne 0)$ then\\
\hspace*{1.3cm} $\bullet$ Write $\tilde{h}=\sum_{(\tilde{g},g)} a_{(\tilde{g},g)}
\cdot \tilde{g}$ with $a_{(\tilde{g},g)} \in \k[x,y]$;\\
\hspace*{1.6cm} (this is possible by Definition~\ref{def:Alg_NF}(2)(3))\\
\hspace*{1.3cm} $\bullet$ $h:=\sum_{(\tilde{g},g)} a_{(\tilde{g},g)} \cdot g$;\\
\hspace*{1.3cm} $\bullet$ $P:=P\cup \{ ( (\tilde{h},h) ,(\tilde{f},f)) \ |
\ (\tilde{f},f) \in S\}$;\\
\hspace*{1.3cm} $\bullet$ $S:=S\cup \{ (\tilde{h}, h) \}$);\\
$\bullet$ Return $S$.
\end{algo}

\begin{rem*}\rm
\begin{itemize}
\item
Notice that if we apply this algorithm to
$G_2=\{0\}$ then we obtain a set of couples $(\tilde{g}, 0)$.
Thus in this situation it is equivalent to Standard.
\item
By construction, for any $(\tilde{g}, g)$
in the output $S$,
$\tilde{g} \in \langle G_2 \rangle \iff g \in \langle G_2 \rangle$.
\item
Applying this algorithm to a basis $G_J$ of
$J$ and a basis $G_\QQ$ of $\QQ$ we get as an
output a set $S$ of couples $(\tilde{g}, g)$. Then
$\tGG=\{\tilde{g} | (\tilde{g}, g) \in S\}$ satisfies
Proposition~\ref{prop:main2} and
$\GG=\{ g | (\tilde{g}, g) \in S\}$
satisfies Corollary~\ref{cor:main2}.
\item
Notice that this algorithm can be used in the general
situation where
we want a standard basis of the sum of two ideals $I_1$ and $I_2$
and such that each $g$ in this basis can be decomposed as $g_1+g_2$
with $g_i \in I_i$.
\end{itemize}
\end{rem*}

Returning to our initial question, we obtain a variant of
PSBmod.

\begin{algo}[$\mathrm{PSBmod}'$]
\label{algo:psbmodpoly}
\rm
\

\noindent
{\rm \bf Input:} $G \subset \k[x,y]$: a finite set,
$G_\QQ \subset \k[y]$: a finite set.\\
{\rm \bf Output:} $\PSBmod'(G,G_\QQ)=(\GG, H)$ with
$\GG$: a pseudo standard basis of $\langle G \rangle$
modulo $\langle G_\QQ \rangle$,
$H \subset \k[y] \smallsetminus \langle G_\QQ \rangle$: a finite set.

\noindent
$\bullet$ Define a global order $\le_0$ on $\N^m$;\\
$\bullet$ Form a block order $\le := (\preceq, \le_0)$;\\
$\bullet$ if $G \subset \langle G_\QQ \rangle$ then
return $(\emptyset, \emptyset)$;\\
$\bullet$ $S:=\mathrm{ModifiedStandard}(G, G_\QQ, \NF)$
where $\NF$ is normal form for $\le$;\\
$\bullet$ $\GG:=\{g \ | \ (\tilde{g},g) \in S, \ g\notin \k[x,y] G_\QQ\}$;\\
$\bullet$ $H:=\emptyset$; for $(g\in \GG)$ do
 $(H:=H \cup \{\lc_\preceq^{\mod\langle G_\QQ \rangle}(g)\})$;\\
$\bullet$ Return $(\GG, H)$.
\end{algo}

To end this part, let us note that if $\preceq$ is not
global one may use a homogenization following Lazard
(see Lemma~\ref{lem:classical}).

\section{Stratification with respect to a constant $\Exp$}

In 4.1, we shall prove Corollary~\ref{cor:stratmain}
and Proposition~\ref{prop:detachable}(2).

In 4.2, 4.3 and 4.4, we propose different variants of
an algorithm illustrating those results.
In 4.2 the algorithm work for a general ring
$\CC$ while the algorithms in 4.3 and 4.4 work when $\CC$
is of the form $\k[y]$.

\subsection{Proof of Cor.~\ref{cor:stratmain}
and Prop.~\ref{prop:detachable}(2)}

Let us recall that we start with an ideal $J\subset \CC[x]$
and $\CC$ is a noetherian integral domain.

We are going to describe a construction by induction on the step $l$.
At each step $l$, we shall construct the following objects:
\begin{itemize}
\item
A finite set $\mathfrak{W}_l$ of triples $(\QQ, h,\GG)$
where $\QQ$ is an ideal of $\CC$, $h\in \CC$ and $\GG$
is finite set in $J$ (the set $\mathfrak{W}_l$ may be empty),
\item
A finite set $\mathfrak{Q}_l$ of ideals $\CC$
(this set may be empty),
\item
An ideal $\II_l$ of $\CC$,
\end{itemize}
with the following properties:
\begin{itemize}
\item[(p1)]
$\spec(\CC)=
(\bigcup_{(\QQ, h, \GG) \in \mathfrak{W}_l} V(\QQ) \smallsetminus V(h))
\cup (\bigcup_{\QQ \in \mathfrak{Q}_l} V(\QQ)) \cup V(\II_l)$,
\item[(p2)]
For any $(\QQ, h, \GG) \in \mathfrak{W}_l$, and for any
specialization $\sigma:\CC \to \K$ such that $\sigma(\QQ)=\{0\}$
and $\sigma(h)\ne 0$,
$\sigma(\GG)$ is $\preceq$-standard basis of
$\K[x]\sigma(J)$,
\item[(p3)]
$J \subset \CC[x] \cdot \II_l$ (i.e. $J$ specializes to zero on
$V(\II_l)$).
\end{itemize}

At step $0$, we set $\mathfrak{W}_0=\emptyset$ and
$\mathfrak{Q}_0=\{(0)\}$ and $\II_0=\langle 1 \rangle$.

Assume the objects of step $l$ are constructed.
If $\mathfrak{Q}_l=\emptyset$ then we stop the
construction. Otherwise we define $\mathfrak{W}_{l+1}$,
$\mathfrak{Q}_{l+1}$ and $\II_{l+1}$ as follows.
Take $\QQ$ in $\mathfrak{Q}_l$.
\begin{itemize}
\item[(A)]
If $J$ is not included in $\CC[x]\cdot \QQ$.\\
Apply Theorem~\ref{theo:main} to $\QQ$.
We obtain $\GG \subset J$ and a finite number of
$h_i \in \CC \smallsetminus \QQ$ ($i=1, \ldots, r$).
We have
$\QQ \subsetneqq \QQ+ \langle h_i \rangle$.\\
Set $h=\prod_{i=}^r h_i$.
Set $\mathfrak{Q}_{l+1}=(\mathfrak{Q}_{l} \smallsetminus \{\QQ\})
\cup \{\QQ + \langle h_i \rangle  | i=1, \ldots, r\}$.\\
Put $\mathfrak{W}_{l+1}=\mathfrak{W}_{l} \cup
\{(\QQ, h, \GG)\}$ and $\II_{l+1}=\II_l$.
\item[(B)]
If $J$ is included in $\CC[x] \QQ$.\\
Set $\mathfrak{Q}_{l+1}:=\mathfrak{Q}_{l} \smallsetminus \{\QQ\}$,
$\mathfrak{W}_{l+1}=\mathfrak{W}_{l}$ and
$\II_{l+1}:=\II_{l} \cap \QQ$.
\end{itemize}
It is clear that at each step $l$, properties (p1), (p2) and
(p3) are satisfied. It is also clear that this construction
is algorithmic if $\CC$ is detachable and intersections are
computable in $\CC$.
Moreover if $\sigma(J)\ne \{0\}$ for any specialization $\sigma$
then $\II_l=\langle 1 \rangle$ for all $l$ (i.e. condition
(B) is never satisfied).
Thus, in order to prove Corollary~\ref{cor:stratmain}
and Proposition~\ref{prop:detachable}(2), it is enough to
prove that there exists $l$ for which 
$\mathfrak{Q}_l$ is empty.

Assume by contradiction that for each $l$, there exists
$\QQ \in \mathfrak{Q}_l$ such that $J \nsubseteq \CC[x]\QQ$.
This will imply the existence of an increasing sequence of
ideals of $\CC$ which contradicts the noetherianity of $\CC$.
Thus there exists $l_0$ such that for all $\QQ \in \mathfrak{Q}_{l_0}$,
$J \subseteq \CC[x] \QQ$. Thus for all steps $l_0, l_0+1, \ldots$,
condition (B) is always satisfied. Thus after a finite number of
steps, $\mathfrak{Q}_l$ becomes empty.

\begin{rem}\

\rm
\begin{enumerate}
\item
Applying this construction, we obtain a union of
$\spec(\CC)$ made of locally closed sets and on each
of these sets $\Exp$ is constant.
Comparing the
values of $\Exp$ on the strata and forming unions of appropriate
strata, we obtain the stratification by a constant $\Exp$.
\item
Notice that given a triple $(\QQ,h,\GG) \in \mathfrak{W}_l$
we may have $V(\QQ)\smallsetminus V(h)=\emptyset$ in $\spec(\CC)$.
Thus, such a triple is useless for the final stratification.
\end{enumerate}
\end{rem}

\subsection{Stratification algorithm 1}

This algorithm consists on a rewriting of the construction
in 4.1.

\begin{algo}[$\mathrm{StratExp1}$]
\label{algo:strat1}
\rm
\

\noindent
{\rm \bf Input:} $G \subset \CC[x]$: a finite set.\\
{\rm \bf Output:} $\mathrm{StratExp}(G)=
(\{(\QQ_1, h_1, \GG_1), \ldots, (\QQ_s, h_s,\GG_s)\}, \II)$;\\
where $\QQ_i \subset \CC$ is an ideal, $h_i \in \CC$,
$\GG_i \subset \CC[x]\cdot G$ is finite
and $\II \subset \CC$ is an ideal.

\noindent
$\bullet$ $\mathfrak{W}:=\emptyset$; $\mathfrak{Q}:=\{(0)\}$;
$\II:=\CC \cdot 1$;\\
$\bullet$ While ($\mathfrak{Q} \ne \emptyset$)\\
\hspace*{.7cm} $\bullet$ Choose $\QQ \in \mathfrak{Q}$;\\
\hspace*{.7cm} $\bullet$ $\mathfrak{Q}:=\mathfrak{Q} \smallsetminus \{\QQ\}$;\\
\hspace*{.7cm} $\bullet$ $(\GG, H):=\PSBmod (G, \QQ)$;\\
\hspace*{.7cm} $\bullet$ if $((\GG, H)\ne (\emptyset, \emptyset))$\\
\hspace*{1.5cm}
$\begin{array}{ll}
\text{then} & \bullet \ h:=\prod_{h' \in H} h';\\
     & \bullet \ \mathfrak{Q}:=\mathfrak{Q} \cup
\bigcup_{h'\in H} \{\QQ+\langle h' \rangle\};\\
     & \bullet \ \mathfrak{W}:=\mathfrak{W} \cup
\{( \QQ, h, \GG  )\}\\
\text{else} & \bullet \ \II:=\II \cap \QQ;
\end{array}$\\
$\bullet$ Return $(\mathfrak{W}, \II)$.
\end{algo}

Here, $(\GG, H)\ne (\emptyset, \emptyset)$
corresponds to condition (A) in 4.1 and
$(\GG, H)= (\emptyset, \emptyset)$ corresponds to condition (B).

\subsection{Stratification algorithm 2}

Here we give a variant of the algorithm above in the case
where $\CC=\k[y]=\k[y_1, \ldots, y_m]$.\\

As we already noticed, in the output of StratExp1 we may
have triples $(\QQ, h, \GG)$ such that
$V(\QQ) \smallsetminus V(h)$ is empty.
In the next variant we may replace the line

$\bullet \ \mathfrak{W}:=\mathfrak{W} \cup
\{( \QQ, h, \GG  )\}$\\
of StratExp1 by

$\bullet \ \text{if } (V(\QQ) \nsubseteq V(h)) \text{ then }
\mathfrak{W}:=\mathfrak{W} \cup \{( \QQ, h, \GG  )\}$\\
the question being how to check the ``if'' condition
in an algorithmic way.

Let us analyse more deeply the construction in 4.1
and show how we may improve it. Let us take the notations
of 4.1. We take $\QQ \in \mathfrak{Q}_l$ such that we
are under condition (A) (i.e. $J \nsubseteq \CC[x] \QQ$).
Applying Theorem~\ref{theo:main}, we obtain $\GG \subset J$
and $h_1, \ldots, h_r \in \CC \smallsetminus \QQ$.
We have:
\[
V(\QQ)=(V(\QQ) \smallsetminus V(h)) \bigsqcup
\bigg( \bigcup_{i=1}^r V(\QQ+\langle h_i \rangle) \bigg) .
\]
The next step consists in adding the triple
$(\QQ,h,\GG)$ to $\mathfrak{W}_l$ and to add
the ideals $\QQ+\langle h_i \rangle $ to
$\mathfrak{Q}_l$.

Although $\QQ \subsetneq \QQ + \langle h_i \rangle$,
we may have $V(\QQ) = V(\QQ +\langle h_i \rangle)$.
Thus in step $l+1$, it would be useless to apply
the construction to $\QQ + \langle h_i \rangle$
if $V(\QQ) = V(\QQ)+\langle h_i \rangle$ i.e.
$V(\QQ) \subseteq V(h_i)$.
Therefore, we would like to replace the line

$\bullet \ \mathfrak{Q}:=\mathfrak{Q} \cup
\bigcup_{h'\in H} \{\QQ+\langle h' \rangle\}$\\
of StratExp1 by

$\bullet \ \mathfrak{Q}:=\mathfrak{Q} \cup
\bigcup_{h'\in H,V(\QQ) \smallsetminus V(h')\ne \emptyset }
\{\QQ+\langle h' \rangle\}$.

Now from an algorithmic point of view how can we check
if a given $h$ is such that $V(\QQ) \smallsetminus
V(h) \ne \emptyset$?
The answer is in the following version of
Hilbert's Nullstellensatz theorem.

\begin{lem}
Let $\k$ be any field. Let $\QQ \subset \k[y]=\k[y_1, \ldots,y_m]$
be an ideal. For any $h\in \k[y]$, we have:
\[
h\in \sqrt{\QQ} \iff V(\QQ) \smallsetminus V(h)
\text{ is empty in } \spec(\k[y]).
\]
\end{lem}
\begin{proof}
We have to prove that $h\in \sqrt{\QQ}$ if and only if
$V(\QQ) \subset V(h)$. The left-right implication is trivial.
Let us assume that $V(\QQ) \subset V(h)$.

Let $\K$ be any algebraically closed field containing $\k$.
Firstly, we have $(\K[y]\QQ)\cap \k[y]=\QQ$. Indeed,
by Remark~\ref{rem:ktoK}, any Gr\"obner basis of $\QQ$
(with respect to a global order)
is a Gr\"obner basis of $\K[y]\QQ$.

Now, given a prime ideal $\PP \subset \K[y]$,
the set $\k[y] \cap \PP$ is a prime ideal of $\k[y]$.
If $\PP$ contains $\K[y]\QQ$ then $\PP \cap \k[y]$
contains $(\K[y]\QQ) \cap \k[y]=\QQ$. Thus the hypothesis
implies that $h \in \PP$. Therefore we have:
$\{y\in \K^m | y\in V_\K(\QQ) \smallsetminus V_\K(h)\}$
is empty. Here $V_\K$ stands for the zero set of.
By the classical Hilbert's Nullstellensatz theorem,
$h^i \in \K[y]\QQ$ for some integer $i$. Finally we obtain:
$h^i \in (\K[y]\QQ) \cap \k[y]=\QQ$.
\end{proof}

Notice that in $\k[y]$, checking whether $h\in \sqrt{\QQ}$
does not require the computation of a Gr\"obner basis
of $\sqrt{\QQ}$, see e.g. \cite[\S 1.8.6]{singular}.

Gathering the previous remarks we obtain the next algorithm.

\begin{algo}[$\mathrm{StratExp2}$]
\label{algo:strat2}
\rm
\

\noindent
{\rm \bf Input:} $G \subset \k[y][x]$: a finite set.\\
{\rm \bf Output:} $\mathrm{StratExp}(G)=
(\{(\QQ_1, h_1, \GG_1), \ldots, (\QQ_s, h_s,\GG_s)\}, \II)$,\\
where $\QQ_i \subset \k[y]$ is a finitely generated ideal,
$h_i \in \k[y]$, $\GG_i \subset \k[y][x] \cdot G$ is finite
and $\II \subset \k[y]$ is an ideal.

\noindent
$\bullet$ $\mathfrak{W}:=\emptyset$; $\mathfrak{Q}:=\{(0)\}$;
$\II:=\k[y] \cdot 1$;\\
$\bullet$ While ($\mathfrak{Q} \ne \emptyset$)\\
\hspace*{.7cm} $\bullet$ Choose $\QQ \in \mathfrak{Q}$; (let $G_\QQ$ denote
a finite basis)\\
\hspace*{.7cm} $\bullet$ $\mathfrak{Q}:=\mathfrak{Q} \smallsetminus \{\QQ\}$;\\
\hspace*{.7cm} $\bullet$ $(\GG, H):=\PSBmod (G, \QQ)$ \text{ or }
$(\GG, H):=\PSBmod' (G, G_\QQ)$;\\
\hspace*{.7cm} $\bullet$ if $((\GG, H)\ne (\emptyset, \emptyset))$\\
\hspace*{1.5cm}
$\begin{array}{ll}
\text{then} & \bullet \ h:=\prod_{h' \in H} h';\\
     & \bullet \ H:=\{h' \in H \ | \ h' \notin \sqrt{Q}\};\\
     & \bullet \ \mathfrak{Q}:=\mathfrak{Q} \cup
\bigcup_{h'\in H} \{\QQ+\langle h' \rangle\};\\
     & \bullet \ \text{if } (h\notin \sqrt{\QQ}) \text{ then }
     \mathfrak{W}:=\mathfrak{W} \cup \{( \QQ, h, \GG  )\}\\
\text{else} & \bullet \ \II:=\II \cap \QQ;
\end{array}$\\
$\bullet$ Return $(\mathfrak{W}, \II)$.
\end{algo}

\subsection{Stratification algorithm 3}

Here we give a usual stratification algorithm for
$\CC=\k[y]$. It uses primary (or prime) decomposition.
In the construction process, all the output tiples
$(\QQ,h,\GG)$ are such that $\QQ$ is prime.
Since $h$ is a product of $h' \in \k[y] \smallsetminus \QQ$,
we shall have $h \notin \QQ$.

We shall give the algorithm without proofs
for correctness and termination since it is
well-known.

\begin{algo}[$\mathrm{StratExp3}$]
\label{algo:strat3}
\rm
\

\noindent
{\rm \bf Input:} $G \subset \k[y][x]$: a finite set.\\
{\rm \bf Output:} $\mathrm{StratExp}(G)=
(\{(\QQ_1, h_1, \GG_1), \ldots, (\QQ_s, h_s,\GG_s)\}, \II)$;\\
where $\QQ_i \subset \k[y]$ is a finitely generated prime ideal,
$h_i \in \k[y]\smallsetminus \QQ$,
$\GG_i \subset \k[y][x]\cdot G$ is finite
and $\II \subset \k[y]$ is an ideal.

\noindent
$\bullet$ $\mathfrak{W}:=\emptyset$; $\mathfrak{Q}:=\{(0)\}$;
$\II:=\k[y] \cdot 1$;\\
$\bullet$ While ($\mathfrak{Q} \ne \emptyset$)\\
\hspace*{.7cm} $\bullet$ Choose $\QQ \in \mathfrak{Q}$; (let $G_\QQ$ denote a finite basis)\\
\hspace*{.7cm} $\bullet$ $\mathfrak{Q}:=\mathfrak{Q} \smallsetminus \{\QQ\}$;\\
\hspace*{.7cm} $\bullet$ $(\GG, H):=\PSBmod (G, \QQ)$ \text{ or }
$(\GG, H):=\PSBmod' (G, G_\QQ)$;\\
\hspace*{.7cm} $\bullet$ if $((\GG, H)\ne (\emptyset, \emptyset))$\\
\hspace*{1.5cm}
$\begin{array}{ll}
\text{then} & \bullet \ h:=\prod_{h' \in H} h';\\
     & \bullet \ \text{Compute prime ideals }
\QQ_1, \ldots, \QQ_r \text{ of } \k[y] \text{ such that}\\
& \hspace{1cm} V(\QQ+\langle h \rangle)=V(\QQ_1)
   \cup \cdots \cup V(\QQ_r);\\
     & \bullet \ \mathfrak{Q}:=\mathfrak{Q} \cup
        \{\QQ_1, \ldots, \QQ_r\};\\
     & \bullet \ \mathfrak{W}:=\mathfrak{W} \cup
\{( \QQ, h, \GG  )\}\\
\text{else} & \bullet \ \II:=\II \cap \QQ;
\end{array}$\\
$\bullet$ Return $(\mathfrak{W}, \II)$.
\end{algo}

\section{Stratification by the local Hilbert-Samuel function}

\begin{proof}[Proof of Corollary \ref{cor:stratHS}]
Recall that we start with a finitely generated ideal
$I \subset \k[x]:=\k[x_1,\ldots,x_n]$ and a field
inclusion $\k\subset \K$.
Let $f_1, \ldots, f_q$ be generators of $I$.
Consider the following ideal
$J=\sum_{l=1}^q \k[x,y]\cdot f_l(x+y)$
where $y$ stands for $(y_1,\ldots,y_n)$.
Take a valuation-compatible order $\preceq$ on the $x^\alpha$'s.
Apply Corollary~\ref{cor:stratmain} to $J\subset \k[y][x]$.
We obtain $\GG_1, \ldots, \GG_r \in J$ and
\[
\spec(\k[y])=(\bigcup_{k=1}^r W_k) \cup V(\II))
\]
where $W_k$ are constructible sets of $\spec(\k[y])$
and $\II \subset \k[y]$ is an ideal.
We have that for any specialization $\sigma$ of $\k[y]$,
if $\sigma(\II)=\{0\}$ then $\sigma(J)=\{0\}$.
Denote by $\overline{\k}$ the algebraic closure of $\k$
and consider the specializations
$\sigma_{y_0}=(\k[y] \to \overline{\k}, P(y) \mapsto P(y_0))$
where $y_0 \in {\overline{\k}}^n$.
If the zeroset $V_{\overline{\k}}(\II) \subset \overline{\k}^n$
is not empty then
for any $y_0\in V_{\overline{\k}}(\II)$ we have $\overline{\k}[x]\cdot
J_{|y=y_0}=\{0\}$ but this is impossible since we implictely
supose $I\ne \{0\}$. Thus $V_{\overline{\k}}(\II)$ is empty, therefore
$1 \in \overline{\k}[y] \II$ i.e. $1\in \II$ and the affine
scheme $V(\II) \subset \spec(\k[y])$ is empty, therefore
$\spec(\k[y])=(\bigcup_{k=1}^r W_k)$.

Now for any specialization $\sigma:\k[y] \to \K$ such that
$\sigma(W_k)=\{0\}$, $\sigma(\GG_k)$ is a $\preceq$-standard
basis of $\K[x] \sigma(J)$ and $\Exp_\preceq(\K[x] \sigma(J))$
does not depend on $\sigma$.
We consider specializations of the form
$\sigma_{x_0}=(\k[y] \to \K, P(y) \mapsto P(x_0))$
where $x_0 \in \K^m$ and use Lemmas~\ref{lem:HS_variable_change}
and \ref{lem:rappelsHS}(2) to conclude.
\end{proof}

\begin{proof}[Proof of Corollary \ref{cor:stratHSgen}]
We sketch the proof since it is similar to the previous one.
Recall that we have an ideal
$I \subset \Z[a,x]$ given by polynomials $f_j=f_j(a,x)$.
Introduce a new set $y$ of indeterminates $y_1, \ldots,y_n$
and consider the ideal
$J=\sum_j \Z[a,y][x]\cdot f_j(a,x+y) \subset \Z[a,y][x]$.
We may apply Corollary~\ref{cor:stratmain} and use the
same arguments as above.
\end{proof}

\section{Examples}

Here we shall give some examples of the computation of a
stratification by the local Hilbert-Samuel function.
These examples (except example 1 treated by hand) were
computed with a program (available on the author's webpage)
written using Risa/Asir computer algebra system \cite{asir}.

In Examples 2 to 5, the output is presented as follows:\\
\noindent
$[[ [\alpha_1, \alpha_2, \ldots ], [q_1(x), q_2(x), \ldots],
  [h_1(x), h_2(x),\ldots]], \ldots ],$

where $\alpha_i \in \N^n$, $q_i, h_i \in \C[x]$.

This means that for
$x_0 \in V(\langle q_1, q_2, \ldots \rangle)
\smallsetminus V(h_1 \cdot h_2 \cdots)$,
the local Hilbert-Samuel function at $x_0$ is equal to that of
the monomial ideal $\langle x^{\alpha_1}, x^{\alpha_2}, \ldots \rangle$.

\subsection{Example 1}

Set $f={x_1}^2+{x_2}^3$ and $I=\C[x_1,x_2] f$.

In this case we shall only use the fact that the Hilbert-Samuel
function associated with $f$ at $x=x_0$ is equal to that
associated with $f(x+x_0)$ at $x=0$.

Let us write $f(x+y)$ as a Taylor series:
\[f(x+y)=(y_1^2 + y_2^3) + (2y_1 x_1 + 3 y_2^2 x_2) +
(x_1^2 + 3 y_2 x_2^2) + (x_2^3).\]
This expansion respects the valuation in $x$. Let us consider a
valuation-compatible order on the $x$-monomials. For $\preceq$, the
leading term of $f(x+y)$ is $1$ and the leading coefficient is
$v=y_1^2 + y_2^3$. On $\C^2 \smallsetminus \{v=0\}$, the
Hilbert-Samuel function is zero.
Now let us work on the space $V:=\{v=0\}$.
Here, working modulo $v$, we can write
\[f(x+y) \equiv (2y_1 x_1 + 3 y_2^2 x_2) +
(x_1^2 + 3 y_2 x_2^2) + (x_2^3).\]
Again we fix a monomial order on $x$ as above. Notice
that we have some freedom: we can choose $\preceq$ in order that
the leading term is $x_1$ or $x_2$ with the corresponding leading
coefficients $2y_1$ or $3 y_2^2$.
Let us choose the leading monomial as $2y_1 \cdot x_1$.
We obtain that on $V \smallsetminus \{(0,0)\}$,
the Hilbert-Samuel function equals that of
$\C[x_1,x_2]/\langle x_1 \rangle$. Finally, it remains $\{(0,0)\}$
(i.e. we work modulo $\langle y_1, y_2 \rangle$)
on which
\[f(x+y) \equiv x_1^2+x_2^3.\]
Finally, we get the stratification
\[\C^2=(\C^2 \smallsetminus V) \cup
(V \smallsetminus \{(0,0)\}) \cup \{(0,0)\}\]
such that one each statum
the local Hilbert-Samuel function is constant
and equal to $\N \ni r\to 0$, $\HSf_{\C[x_1,x_2]/\langle x_1 \rangle }$, and
$\HSf_{\C[x_1,x_2]/\langle x_1^2 \rangle }$ respectively.

\subsection{Example 2}

The same example : $I=\C[x_1,x_2] \cdot (x_1^2+x_2^3)$.
Our program outputs:
\begin{verbatim}
[[[(1)*<<0,0>>],[0],[x1^2+x2^3]],
[[(1)*<<0,1>>,(1)*<<0,1>>],[x1^2+x2^3],[x2,x1]],
[[(1)*<<2,0>>],[x2,x1],[1]]]
\end{verbatim}
We have the following interpretation:
On any point of $V(0)\smallsetminus V(x_1^2+x_2^3)$, the
local Hilbert-Samuel function associated with $I$ is the
same as that of $\C[x_1,x_2]/\langle x_1^0 x_2^0\rangle=
\C[x_1,x_2]/\langle 1 \rangle$. On any point of
$V(x_1^2+x_2^3)\smallsetminus V(x_2 \cdot x_1)$, we get
the Hilbert-Samuel function of $\C[x_1,x_2]/\langle x_2 \rangle$.
On any point of $V(\langle x_2,x_1 \rangle) \smallsetminus V(1)$
(i.e. at $x=(0,0)$), we get the Hilbert-Samuel function of
$\C[x_1,x_2]/\langle x_1^2 \rangle$.
We recover the results of Example 1.

\subsection{Example 3}

Here, we set $f(x_1,x_2,x_3)=x_1^4+x_2^4+x_3 x_1^2 x_2$ and
$I=\C[x_1,x_2,x_3] \cdot f$. The output is~:
\begin{verbatim}
[[[(1)*<<0,0,0>>],[0],[x1^4+x3*x2*x1^2+x2^4]],
[[(1)*<<0,0,1>>,(1)*<<0,0,1>>],[x1^4+x3*x2*x1^2+x2^4],[x1,x2*x1]],
[[(1)*<<2,1,0>>],[x2,x1],[x3]],
[[(1)*<<0,4,0>>],[x3,x2,x1],[1]]]
\end{verbatim}
\noindent
By line 2, we get: On $V(f)\smallsetminus V(x_1 x_2)$, the
local Hilbert-Samuel function is equal to that of
$\C[x_1,x_2,x_3]/\langle x_3 \rangle$.\\
By line 3, we get: On $V(\langle x_1, x_2\rangle)
\smallsetminus V(x_3)$, we have the same Hilbert-Samuel
function as $\C[x_1,x_2,x_3]/\langle x_1^2x_2  \rangle$.\\
By line 3, we get: at $x=(0,0,0)$, the Hilbert-Samuel
function is the same as that of $\C[x_1,x_2,x_3] /\langle
x_2^4 \rangle$.

\subsection{Example 4}

Set $f(x_1,x_2,x_3)=x_1^4+x_2^4+x_3 x_1 x_2$ and $I=\C[x_1,x_2,x_3] \cdot f$.
The program outputs:
\begin{verbatim}
[[[(1)*<<0,0,0>>],[0],[x1^4+x3*x2*x1+x2^4]],
[[(1)*<<0,0,1>>,(1)*<<0,0,1>>],[x1^4+x3*x2*x1+x2^4],[x1,x2*x1]],
[[(1)*<<1,1,0>>],[x2,x1],[x3]],
[[(1)*<<1,1,1>>],[x3,x2,x1],[1]]]
\end{verbatim}

By line 2, we get the Hilbert-Samuel function of
$\C[x_1,x_2,x_3]/\langle x_3 \rangle$ on $V(f) \smallsetminus
V(x_1 x_2)$. By line 3, we get the Hilbert-Samuel function
of $\C[x_1,x_2,x_3]/\langle x_1 x_2\rangle$ on
$V(\langle x_1,x_2 \rangle) \smallsetminus V(x_3)$.
Finally at $x=(0,0,0)$ we get the Hilbert-Samuel function
of $\C[x_1,x_2,x_3]/\langle x_1 x_2 x_3 \rangle$.

\subsection{Example 5}

Here, we set $f_1=x_1-x_2$ and $f_2=x_1 (x_2^2+x_3^3)$
and $I=\C[x_1,x_2,x_3] \{f_1, f_2\}$.
We get the following output (we numbered the lines):
\begin{verbatim}
(1)  [[[(1)*<<0,0,0>>,(1)*<<0,0,0>>],[0],[x1^3+x3^3*x1,-x1+x2]],
(2)  [[(1)*<<0,0,0>>,(1)*<<0,1,0>>],[-x1+x2],[x1^3+x3^3*x1,1]],
(3)  [[(1)*<<0,0,1>>,(1)*<<0,0,1>>,(1)*<<0,1,0>>],
     [-x1+x2,x2^2+x3^3],[x3*x1,x1,1]],
(4)  [[(1)*<<0,1,0>>,(1)*<<3,0,0>>],[x3,x2,x1],[1,1]],
(5)  [[(1)*<<0,1,0>>,(1)*<<1,0,0>>],[x2,x1],[1,x3]],
(6)  [[(1)*<<0,1,0>>,(1)*<<3,0,0>>],[x3,x2,x1],[1,1]],
(7)  [[(1)*<<0,0,0>>,(1)*<<0,0,1>>,(1)*<<0,0,1>>],
     [x1^2+x3^3],[-x1+x2,x3*x1,x1]],
(8)  [[(1)*<<0,0,1>>,(1)*<<0,0,1>>,(1)*<<0,1,0>>],
     [-x1+x2,x2^2+x3^3],[x3*x1,x1,1]],
(9)  [[(1)*<<0,1,0>>,(1)*<<3,0,0>>],[x3,x2,x1],[1,1]],
(10) [[(1)*<<0,0,0>>,(1)*<<3,0,0>>],[x3,x1],[x2,1]],
(11) [[(1)*<<0,1,0>>,(1)*<<3,0,0>>],[x3,x2,x1],[1,1]],
(12) [[(1)*<<0,0,0>>,(1)*<<1,0,0>>],[x1],[x2,x3]],
(13) [[(1)*<<0,1,0>>,(1)*<<1,0,0>>],[x2,x1],[1,x3]],
(14) [[(1)*<<0,1,0>>,(1)*<<3,0,0>>],[x3,x2,x1],[1,1]],
(15) [[(1)*<<0,0,0>>,(1)*<<3,0,0>>],[x3,x1],[x2,1]],
(16) [[(1)*<<0,1,0>>,(1)*<<3,0,0>>],[x3,x2,x1],[1,1]]].
\end{verbatim}

Some lines appear several times (e.g. lines 4, 6, 9, 11, 14
and 16 are equal since they are the termination leaf of several
branches of the tree).
In this result, line 2 (for example) means that on
$V(f_1) \smallsetminus V(f_2)$ the Hilbert-Samuel
function is given by that of $\C[x_1,x_2,x_3]/
\langle 1, x_2 \rangle= \C[x_1,x_2,x_3]/
\langle 1 \rangle$. Line 5 means that on
$V(\langle x_1,x_2\rangle) \smallsetminus \{0\}$ the
Hilbert-Samuel function is given by that
$\C[x_1,x_2,x_3]/ \langle x_1, x_2 \rangle$.

\end{document}